\documentclass[12pt]{article}
\usepackage{latexsym,amssymb,french}

\def\CC{\mathbb C}

\def\PP{\mathbb P}
\def\qed{\hspace*{\fill} $\Box$}
\def\Oo{{\cal O}}
\def\ot{\otimes}
\def\ra{\rightarrow}
\def\11{{\bf 1}}
\def\={:\hspace{-.2mm}=}

\newtheorem{lemma}{Lemme}[section]
\newtheorem{theorem}[lemma]{Th\'eor\`eme}
\newtheorem{definition}[lemma]{D\'efinition}

\newtheorem{proposition}[lemma]{Proposition}

\begin{document}

\title{Le probl\`{e}me de Brill-Noether pour les fibr\'{e}s de
Steiner et application aux courbes gauches}
\author{Ph. Ellia
\footnote{Dipartimento di Matematica, 35 via Machiavelli,
I-44100 Ferrara.}
, A. Hirschowitz
\footnote{Universit\'e de Nice-Sophia Antipolis, 
Parc Valrose, F-06108 Nice Cedex 2.}
, L. Manivel
\footnote {Institut Fourier, 
UMR 5582 du CNRS,
Universit\'e J. Fourier,
BP 74, 
F-38402 Saint-Martin d'H\`eres.}}
\date{}
\maketitle




\section{Introduction} 

On travaille sur $\PP^3=\PP_k^3$ o\`{u} $k$ est un corps 
alg\'{e}briquement clos de caract\'{e}ristique nulle, et l'on 
poursuit (cf. \cite{EH}) l'\'{e}tude des fibr\'{e}s  qui
sont noyaux de morphismes $b.{\Oo}\rightarrow a.{\Oo}(1)$, 
d\'efinis par des matrices de formes lin\'eaires. Ces
fibr\'{e}s, ou plut\^{o}t leurs duaux, ont \'{e}t\'{e} consid\'{e}r\'{e}s,
entre autres, par Dolgachev et Kapranov qui les ont baptis\'{e}s
{\em fibr\'{e}s de Steiner} (cf. \cite{DK}). La raison pour laquelle nous nous
int\'{e}ressons \`{a} ces fibr\'{e}s de Steiner, est qu'ils forment la
classe la plus simple de fibr\'{e}s de syzygies, et que la connaissance de la
cohomologie des fibr\'{e}s de syzygies semble un passage oblig\'{e} pour la
classification des r\'{e}solutions. On aborde ici l'\'{e}tude de la
stratification du champ de ces fibr\'{e}s par la dimension des groupes de
cohomologie (des fibr\'{e}s tordus $E(d)$). L'\'{e}tude de cette
stratification se r\'{e}duit \'{e}videmment \`{a} celle des stratifications
correspondantes sur les espaces de matrices.

Soient $A$ et $B$ des espaces vectoriels de dimension $a$ et $b$.
On note $V$ l'espace $H^0(\PP^3,{\Oo}(1))$, 
et $M$ l'espace des morphismes $m$ de $B$
vers $A\otimes V$ tels que le morphisme correspondant \underline{$m$}$ :b.
{\Oo}\rightarrow a.{\Oo}(1)$ soit surjectif  : $M$ est non vide
d\`{e}s que $b\geq a+3$ (\cite{EH}, Prop.1.1). On note $E_m$ le fibr\'{e}
vectoriel de rang $b-a$, noyau de \underline{$m$}.

Pour $f\geq 1$ on introduit le sous-sch\'{e}ma $W^{f-1}[d]$ de $M$ des
matrices $m$ telles que le corang du morphisme induit  
$$m(d)  : B\otimes S^dV\longrightarrow A\otimes S^{d+1}V$$
soit au moins $f$ (le corang est le minimum des rangs du noyau et du
conoyau). Pour $d=0$, il est facile d'identifier $W^{f-1}[0]$, car
les fibr\'{e}s correspondants admettent un facteur direct trivial de
rang  $f$. Dans ce travail, nous consid\'{e}rons seulement le cas 
suivant, $d=1.$ On fixe donc $a,b$ et $f,$ et 
l'on pose $W\=W^{f-1}[1]$. On fait une fois pour
toutes l'hypoth\`ese $b\geq 5a/2$,
qui correspond au fait que si $m$ est g\'{e}n\'{e}rique dans $M$ alors
$h^0(E_m(1))=4b-10a$ (\cite{EH}, Prop.4.4), et  $m$ appartient \`{a} $W$
si et seulement si $h^1(E_m(1))\geq f$. Dans le m\^eme ordre d'id\'ees,
on se restreint au cas $b \leq 4a$ puisqu'on v\'erifie facilement que tout
fibr\'e de Steiner avec $b > 4a$ est somme directe d'un fibr\'e trivial et
d'un fibr\'e de Steiner avec $b = 4a$.

\smallskip
On peut formuler comme suit le probl\`eme de Brill-Noether pour ces
fibr\'es : 
\begin{enumerate}
\item D\'eterminer l'ensemble $F$ des valeurs  de $f$ pour lesquelles
$W$ est non vide;
\item Pour $f$ dans $F$, \'etudier les composantes irr\'eductibles de $W$;
\item Pour chacune de ces composantes, \'etudier les propri\'et\'es du fibr\'e 
param\'etr\'e par le point g\'en\'erique : en particulier, sa cohomologie.
\end{enumerate}

Dans ce travail, nous \'etudions l'ouvert $W_P$ de $W$ des
matrices $m$ telles que l'image de 
$$m(1)  : B\otimes V\longrightarrow A\otimes S^2V$$
soit de codimension exactement $f$, et o\`u le
morphisme $m \rightarrow Im(m(1))$, \`a valeurs dans la grassmannienne
des sous-espaces de codimension $f$ de $A \otimes S^2V$, est dominant.
Il est assez facile de voir que si $W_P$ est non-vide et non contenu dans
$W^0[0]$, il 
est irr\'eductible.
Notre premier r\'esultat identifie l'ensemble $PF$ des triplets
$(a,b,f) $ (toujours avec $5a/2 \leq b \leq 4a$) 
pour lesquels $W_P$ est non-vide
et non contenu dans $W^0[0]$, et assure que, pour $(a,b,f)$ dans $PF$, le
fibr\'e param\'etr\'e par le point g\'en\'erique de $W_P$ a la cohomologie
attendue:

\begin{theorem}
Soient $a$ et $b$ deux entiers positifs avec
$5a/2 \leq b \leq 4a$, et soit $f$ positif. Pour que $W_P$ ne soit pas contenu
dans $W^0[0]$, il faut et il suffit que $f$ soit inf\'erieur \`a
$a - b/4$. Dans ce cas, $W_P$ est irr\'eductible
et le fibr\'e $E_\mu$ param\'etr\'e par son point
g\'en\'erique $\mu$ a la cohomologie attendue en ce sens 
que, pour $k \neq 1$, au plus un des quatre
nombres $h^i(E_\mu(k))$ est non nul.
\end{theorem}

Notre m\'{e}thode de d\'emonstration
consiste, comme dans \cite {EH},
 \`{a} consid\'{e}rer les matrices $m$ telles que 
l'image de $m(1)$ soit contenue dans un sous-espace de codimension donn\'ee
$f$ de $A\otimes S^2V$, et \`{a} majorer le lieu des mauvaises matrices en le
stratifiant judicieusement.

Le r\'esultat pr\'ec\'edent est insuffisant pour l'application aux courbes
qui est la motivation initiale de ce travail. Il nous faut en effet
savoir, pour le fibr\'e $E_\mu$ du th\'eor\`eme pr\'ec\'edent, si
le morphisme d'\'evaluation $H^0(E_\mu(1))\otimes \Oo \ra E_\mu(1)$
est {\em versel}, i.e. si les lieux o\`u son rang est
constant sont lisses de la codimension attendue(\cite {EH}). Ce probl\`eme 
prend deux tournures assez diff\'erentes selon que la dimension
$4b - 10a +f$ de $H^0(E_\mu(1))$ est inf\'erieure ou sup\'erieure
au rang $b-a$ de $E_\mu(1)$, et nous n'abordons ici que le second cas
(l'autre devrait permettre de construire des courbes gauches \`a {\em monade}
lin\'eaire). On suppose donc $f \geq 9a - 3b$ et on obtient le r\'esultat
suivant.

\begin{theorem} \label {versel}
On suppose $5a\le 2b$, et 
$$9a-3b\le f\le \min(\frac{3a-1}{11}, \frac{13a-4b-5}{5},
\frac{16a-5b-5}{2}).$$
Alors,  si $\mu $ d\'{e}signe le point g\'{e}n\'{e}rique de $W_P$,
le fibr\'{e} $E_\mu (1)$ est d'\'{e}valuation verselle  : 
les lieux o\`{u} le morphisme
d'\'{e}valuation $H^0(E_\mu (1))\otimes {\Oo}\rightarrow E_\mu (1)$
est de rang constant sont lisses, de la codimension attendue.
\end{theorem}

Ce th\'eor\`eme se d\'emontre en poussant sensiblement plus loin la
m\'ethode utilis\'ee pr\'ec\'edemment. Les bornes que nous imposons
\`a $f$ ne semblent pas trop restrictives puisque, dans notre application,
nous obtenons pratiquement le r\'esultat optimal.

L'application mentionn\'{e}e dans le titre concerne la classification des
courbes (g\'{e}n\'{e}riques) de $\PP^3$, et plus pr\'{e}cisement
l'existence de courbes g\'{e}n\'{e}riques ''principales''. Il est bien connu
qu'il existe un domaine de $(d,g)$ (''domaine A'') o\`{u} l'on
s'attend \`a
l'existence d'une unique composante irr\'{e}ductible du sch\'{e}ma de
Hilbert $H(d,g)$ (composante ''principale''), telle que la courbe
g\'{e}n\'{e}rique correspondante ait toutes les ''bonnes''
propri\'{e}t\'{e}s voulues. En particulier la r\'{e}solution libre minimale
de l'id\'{e}al d'une telle courbe devrait \^{e}tre lin\'{e}aire ou presque
lin\'{e}aire (il n'y a que quatre formes possibles). Une courbe
g\'{e}n\'{e}rique $C$ de $\PP^3$ est dite \`{a} r\'{e}solution
lin\'{e}aire si son id\'{e}al ${\cal I}_C$ admet une r\'{e}solution de la
forme :
$$0\longrightarrow a.{\Oo}(-s-2)\longrightarrow b.{\Oo}(-s-1)\longrightarrow
{\Oo}(-s)\longrightarrow {\cal I}_C\longrightarrow 0.$$

Observons que si le morphisme $a.{\Oo}(-s-2)\rightarrow b.{\Oo}
(-s-1)$ est vu comme morphisme $m$ de $B$ vers $A\otimes V,$ alors $m(1)$
est de rang au plus $3b+a-1$ (car le dual du fibr\'{e} des syzygies tordu
par $-s$ a au moins $b-a+1$ sections). En particulier si $b\leq 3a$, alors 
$m$ ne peut \^{e}tre g\'{e}n\'{e}rique.

Nous dirons que la courbe g\'{e}n\'{e}rique 
$C$ ci-dessus est \`{a} {\em r\'{e}solution
lin\'{e}aire pr\'{e}dominante} (rlp), de type $(a,b),$ si l'image de $m(1)$
est le sous-espace g\'{e}n\'{e}rique (en ce sens que son image dans le
quotient de la grassmannienne par $Gl(A)$ est le point g\'{e}n\'{e}rique) de
dimension $\min (10a,3b+a-1)$ de $A\otimes S^2V$. Rappelons que les courbes
lin\'{e}aires dominantes introduites dans \cite{EH} pour $b>3a$ ont $m(1)$
surjectif et donc sont rlp. En utilisant les th\'{e}or\`{e}mes
pr\'{e}c\'{e}dents avec $f=9a-3b+1$ nous obtenons  :

\begin{theorem}
Pour $b>(3-\frac{1}{11})a+\frac{9}{11}$ et $a\ge 7$, 
 il existe dans $\PP^3$ une
unique courbe g\'{e}n\'{e}rique lisse, g\'{e}om\'{e}triquement connexe, 
\`a r\'esolution lin\'eaire pr\'edominante de type $(a,b).$
\end{theorem}

Cet \'enonc\'e est tout pr\`es d'\^etre optimal 
puisque la condition $b>(3-\frac{1}{11})a+\frac{1}{3}$
est n\'ecessaire (remarque 7.3). On observera 
par ailleurs que la construction des pr\'{e}sentes courbes est sensiblement
plus d\'{e}licate que celle des courbes \`{a} r\'{e}solution lin\'{e}aire
dominante de \cite{EH}. Ceci correspond au fait que nos nouvelles courbes
sont construites comme lieux de d\'{e}pendance de sections de fibr\'{e}s
qui, contrairement au cas de \cite{EH}, {\it ne sont pas engendr\'{e}s
par leurs sections} (remarque 4.3). 

Remarquons \'egalement que l'on ne peut esp\'{e}rer que toutes les
courbes g\'{e}n\'{e}riques principales \`{a} r\'{e}solution lin\'{e}aire
soient rlp  : en effet, les composantes rlp du sch\'{e}ma de Hilbert sont
unirationnelles, tandis que les composantes principales du domaine de
Brill-Noether ne le sont en g\'{e}n\'{e}ral pas (\cite{HM}, \cite{H}).

Enfin, nous pensons que les techniques developp\'{e}es dans
ce travail (ainsi que dans \cite{EH}) devraient s'appliquer aux cas
''mixtes'' (r\'{e}solutions ''quasi-lin\'{e}aires''), dans des domaines
analogues, et devraient donc permettre de r\'{e}soudre le probl\`{e}me de la
composante principale dans une frange enti\`{e}re (''au sommet'') 
du domaine A.


\section{ F-formes et transport d'\'equations}
\subsection{$F$-formes trilin\'{e}aires partiellement sym\'{e}triques}

Soient $A,F$ et $V$ trois espaces vectoriels de dimension $a,f$ et $n+1$
respectivement. On appelle $F$-forme tout morphisme \`{a} valeur dans $F$.
On dira qu'une $F$-forme trilin\'{e}aire $f$ sur $A\times V\times V$ est
{\em partiellement sym\'{e}trique} si elle v\'{e}rifie l'identit\'{e} 
$f(x,y,z)=f(x,z,y)$. Une telle $F$-forme trilin\'{e}aire donne naissance
\`{a} des $F$-formes lin\'{e}aires sur $A\otimes S^2V$ et
$A\otimes V\otimes V$, la deuxi\`{e}me se d\'{e}duisant de la
premi\`{e}re par composition avec le produit $A\otimes V\otimes
V\longrightarrow A\otimes S^2V$. Elle donne aussi naissance \`{a} une
application lin\'{e}aire de $A\otimes V$ vers $V^{\vee}\otimes F$,
associ\'{e}e \`{a} l'application bilin\'{e}aire de $A\times V$ vers 
$V^{\vee}\otimes F={\rm Hom}\, (V,F) :(x,y)
\longrightarrow (z\longrightarrow f(x,y,z))$.
Inversement toute $F$-forme lin\'{e}aire sur $A\otimes S^2V$ provient
d'une $F$-forme trilin\'{e}aire partiellement sym\'{e}trique, qu'on obtient
en composant avec le produit naturel de $A\times V\times V$ vers 
$A\otimes S^2V$.

Dans la suite on rencontrera d'abord des $F$-formes lin\'{e}aires sur $
A\otimes S^2V$. Si $g$ est une telle $F$-forme lin\'{e}aire sur $
A\otimes S^2V$, on notera $g^{*}$ le morphisme correspondant de $
A\otimes V$ vers $V^{\vee}\otimes F$ et $g^{\prime }$ celui de $
F^{\vee}\otimes S^2V$ vers $A^{\vee}$.
Les morphismes $g$ et $g^{*}$ sont li\'{e}s par le diagramme commutatif :
$$ (D.1)\hspace{1cm}
\begin{array}{cccc}
& A\otimes V & = & A\otimes V \\
& \hspace*{4mm}\downarrow {\scriptstyle j^*} & 
 & \hspace*{4mm}\downarrow {\scriptstyle g^*}  \\
& A\otimes S^2V\otimes V^{\vee} & 
\stackrel{{\scriptstyle g\otimes 1}}{\longrightarrow}
& V^{\vee}\otimes F \end{array}$$
o\`{u} $j$ est le compos\'{e} de : $\11\otimes \pi \otimes \11 :A\otimes
(V\otimes V)\otimes V^{\vee}\longrightarrow A\otimes S^2V\otimes V$
avec $\11\otimes j  : (A\otimes V)\otimes k\longrightarrow A\otimes
V\otimes V\otimes V^{\vee}$.
Le diagramme pr\'{e}c\'{e}dent exprime $g^{*}$ en fonction de $g$. Pour
exprimer inversement $g$ en fonction de $g^{*}$, on dispose du diagramme
commutatif suivant (qui est une cons\'{e}quence imm\'{e}diate du cas
classique avec $A=F=k$) :
$$ (D.2)\hspace{1cm}
\begin{array}{cccc}
& A\otimes V\otimes V & 
\stackrel{g^*\otimes 1}{\longrightarrow}
 & F\otimes V^{\vee}\otimes V \\
& {\scriptstyle 1\otimes s}\downarrow\hspace*{6mm} &  & 
\hspace*{6mm}\downarrow {\scriptstyle 1\otimes t} \\
& A\otimes S^2V & \stackrel{g}{\longrightarrow} & F\otimes k
\end{array}$$
o\`{u} $s$ et $t$ sont les morphismes naturels de produit et de trace.

\subsection{$F$-formes trilin\'{e}aires tr\`{e}s partiellement
sym\'{e}triques}

On g\'{e}n\'{e}ralise un peu ce qui pr\'{e}c\`{e}de. Si
$H$ est un hyperplan de $V$, 
on dit qu'une $F$-forme trilin\'{e}aire sur $A\times
H\times V$ est {\em tr\`{e}s partiellement sym\'{e}trique} 
(tps) si sa restriction
\`{a} $A\times H\times H$ est partiellement sym\'{e}trique. Une $F$-forme
tps donne naissance \`{a} une $F$-forme lin\'{e}aire sur $A\otimes
H\otimes V$ qui se factorise \`{a} travers $A\otimes H.V$, et
\`{a} un morphisme de $A\otimes V$ vers $H^{\vee}\otimes V$. Inversement
toute $F$-forme lin\'{e}aire $f$ sur $A\otimes H.V$ provient d'une 
$F$-forme tps, et on note $f^{*}$ le morphisme associ\'{e} de $A\otimes V$
vers $H^{\vee}\otimes F$. Comme plus haut $f$ et $f^{*}$ sont li\'{e}s par les
diagrammes commutatifs  :
$$ (D.3)\hspace{1cm}
\begin{array}{cccc}
& A\otimes V & = & A\otimes V \\
& \hspace*{5mm}\downarrow {\scriptstyle j^*} &  & 
\hspace*{5mm}\downarrow {\scriptstyle f^*} \\
& A\otimes H.V\otimes H^{\vee} & \stackrel{f\otimes 1}{\longrightarrow} &
H^{\vee}\otimes F \end{array}$$
o\`{u} $j^{*}$ est le compos\'{e} de
$\11\ot\pi\ot\11  : A\otimes (V\otimes H)\otimes
H^{\vee}\rightarrow A\otimes H.V\otimes H^{\vee} $
avec $\11\otimes j :(A\otimes V)\otimes k\rightarrow A\otimes
V\otimes H\otimes H^{\vee},$ et 
$$ (D.4)\hspace{1cm}
\begin{array}{cccc}
& A\otimes V\otimes H & \stackrel{f^*\otimes\11}{\longrightarrow} & 
F\otimes H^{\vee}\otimes H \\
& \hspace*{6mm}\downarrow {\scriptstyle \11\otimes s} &  & 
\hspace*{6mm}\downarrow {\scriptstyle \11\otimes t} \\
& A\otimes H.V & \stackrel{f}{\longrightarrow} & F\otimes k
\end{array} $$

Les $F$-formes constituent un moyen de calculer autour des objets qui nous
int\'{e}ressent, \`{a} savoir les sous-espaces de $A\otimes S^2V$. Si $Z$
est un tel sous-espace, on note $F_Z$ le quotient correspondant, $Z^{*}$ le
noyau du morphisme correspondant de $A\otimes V$ vers $F_Z\otimes V^{\vee}$.

\subsection{Transports d'\'{e}quations}

Dans ce paragraphe on montre qu'une \'{e}quation, portant sur $m(1)$ (ou sur
$m_H(1)$), de la forme $g\circ m(1)=0$, \'{e}quivaut \`{a} un syst\`{e}me
d'\'{e}quations de la forme $g^{*}\circ m=0$. L'avantage est que le rang de
ce nouveau syst\`{e}me est plus accessible (on l'appelle 
le $V^{*}-$rang de $g$, cf. d\'efinition 3.1).
        
\begin{proposition}

Soit $g$ une $F$-forme lin\'{e}aire sur $A\otimes S^2V$. Le morphisme $
m  : B\rightarrow A\otimes V$ v\'{e}rifie $g\circ m(1)=0$ si et seulement
si $g^{*}\circ m=0$.

\end{proposition}

\noindent {\em D\'{e}monstration :} 
En utilisant (D.1), on obtient le diagramme commutatif :
$$\begin{array}{ccccc}
B & \stackrel{m}{\longrightarrow} & A\otimes V &  &  \\
\hspace*{6mm}\downarrow {\scriptstyle \11\otimes j} &  & 
{\scriptstyle j^*} \downarrow\hspace*{5mm} & \searrow^{g^{*}} &  \\
B\otimes V\otimes V^{\vee} & \stackrel{m(1)\otimes 1}{\longrightarrow} & 
A\otimes S^2V\otimes V^{\vee} & \stackrel{g\otimes 1}{\longrightarrow} & 
V^{\vee}\otimes F. \end{array}$$
L'hypoth\`{e}se $g\circ m(1)=0$ implique que $(g\otimes 1)\circ (m(1)\otimes
1)$ est nul. Par suite $g^{*}\circ m$ est nul. Inversement on construit
\`{a} partir de (D.2) le diagramme commutatif suivant  :
$$\begin{array}{ccccc}
B\otimes V & \stackrel{m\otimes 1}{\longrightarrow} & A\otimes V\otimes V &
\stackrel{g^{*}\otimes 1}{\longrightarrow} & F\otimes V^{\vee}\otimes V \\
& {\scriptstyle m(1)}\searrow\hspace*{6mm} & \hspace*{6mm}\downarrow 
{\scriptstyle \11\otimes s} &  & \hspace*{6mm}\downarrow 
{\scriptstyle \11\otimes t} \\
&  & A\otimes S^2V & \stackrel{g}{\longrightarrow} & F\otimes k
\end{array}$$
qui montre que si $g^{*}\circ m$ est nul, alors $g\circ m(1)$ l'est aussi.
\qed

\begin{proposition}
Soit $f$ une forme lin\'{e}aire sur $A\otimes H.V$. Alors le morphisme $m$
v\'{e}rifie $f\circ m_H(1)=0$ si et seulement si $f^{*}\circ m=0$.
\end{proposition}

\noindent {\em D\'{e}monstration :} Analogue \`{a} la pr\'{e}c\'{e}dente, en
utilisant (D.3) et (D.4) au lieu de (D.1) et (D.2).\qed

\begin{proposition}

Soit $Z$ un sous-espace de codimension $f$ de $A\otimes S^2V$ et 
$Z^{\prime}$ l'intersection, suppos\'{e}e transverse, de $Z$ avec $A\otimes
H.V$. Soit $T$ un sous-espace de $A\otimes H.V$ contenu dans
$Z^{\prime }$. On pose 
$F^{\prime } \= (A\otimes H.V)/T$ et $F \=(A\otimes S^2V)/Z$, et l'on
consid\`ere le morphisme correspondant
$$f_{Z,T}^{*}  : A\otimes V\longrightarrow
(V^{\vee}\otimes F)\oplus (H^{\vee}\otimes F^{\prime }).$$ 
Alors, pour que l'image de $m(1)$ soit contenue dans $Z$ et celle de 
$m_H(1)$ dans $T$, il faut et il suffit que $f_{Z,T}^{*}\circ m=0$.

\end{proposition}

\noindent {\em D\'{e}monstration :} C'est une cons\'{e}quence 
imm\'{e}diate des deux propositions pr\'{e}c\'{e}dentes. \qed

\section{Rangs auxiliaires et stratifications}

Pour comprendre les sous-espaces lin\'eaires de $A\otimes S^2V$, nous 
aurons besoin d'introduire des quantit\'es auxiliaires, qui seront les
rangs de certaines applications associ\'ees. 
  
\begin{definition}
(du $V^{*}$-rang). Si $g$ est une $F$-forme sur $A\otimes S^2V$, son 
$V^{*}$-rang est par d\'{e}finition le rang de $g^{*} :A\otimes V\rightarrow
F\otimes V^{\vee}$.
Si $Z$ est un sous-espace lin\'{e}aire de $A\otimes S^2V$ on d\'{e}finit son
$V^{*}$-rang comme le $V^{*}$-rang de son conoyau.
\end{definition}

\begin{lemma}
Si $a\geq f$, le sous-espace lin\'{e}aire g\'{e}n\'{e}rique de
codimension $f$ de $A\otimes S^2V$ est de $V^{*}$-rang maximum,
\`a savoir $4f$.
\end{lemma}

\noindent {\em D\'{e}monstration :} Soit $F$ un espace vectoriel 
de dimension $f$, et $x :A\rightarrow F$ surjectif d'une part,  
$y :V\rightarrow V^{\vee}$ sym\'{e}trique non-d\'{e}g\'{e}n\'{e}r\'{e}e
d'autre part. Alors $x\otimes y$ provient d'une $F$-forme
lin\'{e}aire sur $A\otimes S^2V$, dont le noyau est de codimension $f$
et de $V^{*}$-rang maximum, c'est-\`a-dire $4f$.\qed

\begin{definition}
(du $Z$-rang). Si $T$ et $Z$ sont deux sous-espaces lin\'{e}aires de $A\otimes
S^2V,$ avec $T\subset Z$, on appelle $Z$-rang de $T$ le $V^{*}$-rang de $T$
diminu\'{e} du $V^{*}$-rang de $Z$.
\end{definition}

On va maintenant  stratifier par le $Z$-rang
l'espace projectif des hyperplans du sous-espace lin\'eaire g\'en\'erique
$Z$ de codimension $f$ de $A\otimes S^2V$. On notera 
donc $G_r$ le sous-sch\'ema des hyperplans de $Z$ de $Z$-rang au plus $r$.
Dans le reste de cette section, on majore la dimension de $G_r$.

\begin{proposition} 
Supposons  $2a>5f$. Alors $G_0$ est vide, et on a
$$\begin{array}{lcl}
{\rm codim}\; G_1 & \geq & \min (12a-12f-3,11a-5f-5,9a-f-3), \\
{\rm codim}\; G_2 & \geq & \min(7a-f-4,11a-7f-4), \\
{\rm codim}\; G_3 & \geq & 4a-4f-3.
\end{array}$$
\end{proposition} 

\noindent {\em D\'emonstration.}
Soit $\Phi $ un morphisme surjectif de $A\otimes S^2V$ vers $k^f$ dont $Z$ est 
le noyau. Au lieu de $G_r$, il nous
suffit de stratifier le sch\'{e}ma, de dimension $10a$, des (noyaux des)
morphismes surjectifs $g$ de $A\otimes S^2V$ vers $k^{f+1}$ de la forme 
$(\Phi ,-)$. Choisissons des bases $(e_j)$ sur $A$ et $(v_p)$ sur $V$.
La matrice sym\'{e}trique $X$ de $\Phi $ s'\'{e}crit $x_{pqjs}$ avec $1\leq
p,q\leq 4$, $1\leq j\leq a$ et $1\leq s\leq f$. Celle de $g$ s'\'{e}crit 
$g_{pqjs}$ avec $1\leq p,q\leq 4$, $1\leq j\leq a$ et $1\leq s\leq
f+1$, et pour $s\leq f$, on a $g_{pqjs}=x_{pqjs}$. 
On pose $g_{pqj}=g_{pqj,f+1}$. Le $Z$-rang de $
\ker (g)$ est donc le surcro\^{\i}t de rang qu'apportent les ''lignes'' $
g_{p..} $ aux ''lignes'' $x_{p..s}$ de $\Phi $. Comme dans \cite{EH}, les
majorations de dimension se font localement.

\medskip\noindent {\bf Cas $r=0$}. 
On majore la dimension de l'ensemble $\Gamma$ 
des couples $(X,C)$ tels que $g_{pqj}= \sum_{rs}c_{prs}x_{qrjs}$ 
soit sym\'etrique. Si l'on projette sur $C$, alors $A$ se factorise, 
c'est-\`a-dire que la dimension des fibres est lin\'eaire en celle 
de $A$, de sorte qu'il suffit de la calculer lorsque $A$
est de dimension $1$. 

On traite d'abord le cas  $f=1$. La dimension
de l'ensemble des $X$ convenables pour $C$ fix\'e ne d\'epend que de la 
forme de Jordan de cette matrice $4\times 4$  : en effet, matriciellement,
$G=CX$ est sym\'etrique si et seulement si 
$PGP^t= (PCP^{-1})(PXP^t)$ l'est aussi.
On stratifie donc selon la forme de Jordan $J(C)$ de $C$, cod\'ee par
la taille des blocs de Jordan de chaque valeur propre. Le tableau 
suivant donne les dimensions $O(C)$  et $S(C)$ des strates  et des espaces 
de solutions correspondants (les barres verticales
s\'eparent les tailles des blocs 
correspondant \`a des valeurs propres distinctes).
\medskip\begin{center}
\begin{tabular}{||c|c|c||c|c|c||} \hline
$J$ & $O$ & $S$ &  $J$ & $O$ & $S$ \\ \hline
$1|1|1|1$ & 16 & 4  & $21|1$ & 12 & 5 \\  
$2|1|1$ & 14 & 4  & $31$ & 11 & 5 \\ 
$2|2$ & 14 & 4 & $11|11$ & 10 & 6 \\
$3|1$ & 14 & 4  & $22$ & 9 & 7 \\
$4$ & 13 & 4 & $111|1$ & 8 & 7 \\
$11|1|1$ & 13 & 5 & $211$ & 7 & 7 \\
$2|11$ & 12 & 5 & $1111$ & 1 & 10 \\ \hline
\end{tabular}
\end{center}\medskip
La v\'erification de ce tableau est \'el\'ementaire, et consiste \`a
 constater que, dans chaque cas, les composantes non identiquement 
nulles de la partie antisym\'etrique de $CX$ forment un ensemble
de $10-S$ formes lin\'eaires ind\'ependantes.

Lorsque $f$ est quelconque, on consid\`ere $c$ comme une suite
$C= (C^1,\ldots ,C^f)$ de $f$ matrices $4\times 4$, et l'on 
stratifie selon le maximum $O(C)$ des $O(C^i)$, tandis
qu'on note $S(C)$ le minimum des $S(C^i)$.
Les dimensions $O_f(C)$ et $S_f(C)$ des strates 
et des espaces de solutions 
correspondants v\'erifient \'evidemment
$$O_f=  f\times O,\quad S_f\leq S+10(f-1).$$
On observe maintenant que les fibres de la restriction \`a $\Gamma$ de
la premi\`ere projection sont de dimension au moins $f$,
\`a cause de l'invariance par
$c_{pqs}\rightarrow c_{pqs}+\gamma_s\delta _{pq}$.
On montre alors que la seule strate de $\Gamma$ dans laquelle $X$ est
g\'en\'erique est la derni\`ere, celle o\`u $O=1$. En effet, pour chacune
des autres strates, on a $(10-S)a>(O-1)f$
gr\^ace \`a l'hypoith\`ese sur $f$, et donc $10af>aS_f+O_f-f.$
Quant \`a la derni\`ere strate, elle ne nous concerne pas puisque ses
points ne d\'efinissent pas un sous-espace de codimension $f+1$ de
$A\otimes S^2V$. On a donc montr\'e que $G_0$ est vide.

\medskip\noindent {\bf Cas $r=1$}. 
Quitte \`a changer de base, on peut supposer ici que 
le surcro\^{\i}t de rang est d\^u \`a $g_{1..}$,
d'o\`u les relations de d\'ependance, pour $2\leq p\leq 4$,
$$g_{p..}= \alpha _pg_{1..}+\sum_{qs}c_{pqs}x_{q..s}.$$
Les relations de sym\'etrie $g_{pq.}=g_{qp.}$ pour $2\leq p<q\leq 4$
s'\'ecrivent
$$\alpha _p\sum_{rs}c_{qrs}x_{1r.s}+\sum_{rs}c_{prs}x_{qr.s}=  
\alpha _q\sum_{rs}c_{prs}x_{1r.s}+\sum_{rs}c_{qrs}x_{pr.s}.$$
On estime alors la dimension du sch\'ema $\Gamma_1$ des
$(\alpha, c , g, x)$ v\'erifiant les \'equations ci-dessus en stratifiant
le sch\'ema des  $(\alpha ,c)$ selon la dimension des solutions en $x$  : 
$A$ se factorise. 

G\'en\'eriquement, le syst\`eme pr\'ec\'edent est de 
rang $3$, et la strate correspondante de $\Gamma_1$ est de dimension
$3$ (pour les $\alpha$) plus $12f$ (pour $c$) plus $a$ (pour $g_{11.}$)
plus $a(10f-3)$ (pour les $x$), et donc son image est de codimension
au moins
$12a-12f-3$ dans la vari\'et\'e des morphismes de $A \otimes S^2V$
vers  $k^{f+1}$. Par suite la strate correspondante de $G_1$ est de 
codimension au moins $12a-12f-3$. 

Si le syst\`eme est de rang deux, 
la strate correspondante de $\Gamma_1$ est de dimension $3$ 
(pour les $\alpha$) plus $2$ (pour la relation de d\'ependance entre 
les trois \'equations) plus  $5f$ (pour $c$, en effet, on v\'erifie
facilement que selon la relation de d\'ependance, il y a moyen de
choisir $7f$ parmi les $c_{pqs}$ qui se calculent en fonction des
autres, des $\alpha$ et de la relation) plus $a$ (pour $g_{11.}$)
plus $a(10f-2)$ (pour les $x$).
La dimension de la strate correspondante est donc major\'ee par 
$5f-a+5$. 

On remarque enfin que le syst\`eme n'est jamais de rang un,
et que le rang z\'ero ne laisse subsister que $a+f+3$ param\`etres (\`a
savoir les $g_{11.}$, les $c_{11.}$ et les $\alpha$),
d'o\`u une derni\`ere strate de dimension $a+f+3$.

\medskip\noindent {\bf Cas $r=2$}.
On raisonne comme pr\'ec\'edemment, en supposant
le surcro\^{\i}t de rang d\^u \`a $g_{1..}$ et $g_{2..}$,
les relations de d\'ependance \'etant
$$\begin{array}{rcl}
g_{3..} & =  & \alpha _{31}g_{1..}+\alpha _{32}g_{2..}+\sum_{rs}c_{3rs}
x_{r..s}, \\
g_{4..} & =  & \alpha _{41}g_{1..}+\alpha _{42}g_{2..}+\sum_{rs}c_{4rs}
x_{r..s}.
\end{array}$$
D'o\`u l'unique relation de sym\'etrie
$$\begin{array}{l}
\alpha _{31}\sum_{rs}c_{4rs}x_{r1.s}+\alpha _{32}\sum_{rs}c_{4rs}x_{r2.s}
+\sum_{rs}c_{3rs}x_{r4.s} \\
\hspace*{3cm} =  \alpha _{41}\sum_{rs}c_{3rs}x_{r1.s}
+\alpha _{42}\sum_{rs}c_{3rs}x_{r2.s}+\sum_{rs}c_{4rs}x_{r3.s}.
\end{array}$$
Dans le cas g\'en\'erique, cette relation est non triviale,
 et l'on trouve une strate de $\Gamma_2$ de dimension 
$4$ (pour les $\alpha$) plus $8f$ (pour $c$) plus $3a$ (pour $g_{11.},
g_{12.}, g_{22.}$)
plus $a(10f-1)$ (pour les $x$), dont la trace sur
la vari\'et\'e des morphismes prolongeant $f$ 
de $A \otimes S^2V$ vers $k^{f+1}$ est donc  de dimension au plus
$2a+8f+4$, c'est-\`a-dire de codimension au moins $8a-8f-4$.
Par suite la strate correspondante de $G_2$ est \'egalement de codimension au
moins $8a-8f-4$.

Quand l'\'equation en $x$ ci-dessus est identiquement nulle, on a
$$c_{34s}= c_{43s}= 0,\quad c_{33s}= c_{44s} \= c_s, \quad
{\rm et}\;\; c_{pqs}= -\alpha _{pq}c_s \;\; {\rm si}\;\;q \leq 2<p.$$
On trouve ainsi une strate de $\Gamma_2$ de dimension $3a+f+4+10af$, 
dont la trace sur la vari\'et\'e des morphismes prolongeant $f$ 
de $A\otimes S^2V$ vers $k^{f+1}$ est donc  de dimension au plus
$3a+f+4$, c'est-\`a-dire de codimension au moins $7a-f-4$.
Par suite la strate correspondante de $G_2$ est \'egalement de 
codimension au moins $7a-f-4$.

\medskip\noindent {\bf Cas $r=3$.} 
 Quitte \`{a} changer de base dans $V$, on peut supposer que le
surcro\^{\i}t de rang est d\^{u} \`{a} $g_{1..}, g_{2..},g_{3..}$. 
L'ouvert correspondant est param\'{e}tr\'{e} par les $g_{pq.}$ pour 
$1\le p\le q\le 3$, et les coefficients exprimant $g_{4..}$ en fonction
lin\'{e}aire des $x_{p..s}$ et de $g_{1..}, g_{2..},g_{3..}$.
On majore ainsi la dimension par $6a+4f+5$.\qed

\section{La composante PW}

Pour $m$ dans ${\cal M}={\rm Hom}\;(B,A\otimes V)$, 
on d\'{e}signera par $Z_m$
l'image de $m(1)$ et par $f_m :A\otimes S^2V\rightarrow F_m$ son conoyau. En
outre, $Z_m^{*}$ d\'{e}signera le noyau du morphisme correspondant de 
$A\otimes V$ vers $F_m\otimes V^{\vee}$ (cf. 2.2).

Soit ${\cal W} \={\cal W}^{f-1}[1]$ le
sous-sch\'{e}ma de ${\cal M}$ des matrices telles que le corang (le
minimum des rangs du noyau et du conoyau) de $m(1)$ soit au moins $f$. 
Si $4b\geq 10a,$ ${\cal W}$ est le sous-sch\'{e}ma des $m$
dans ${\cal M}$ tels que ${\rm codim}\; Z_m\geq f$. 
On note $W$ la trace de ${\cal W}$ sur l'ouvert $M$ de ${\cal M}$.

\begin{proposition}
On suppose que $2b\ge 5a$ et $2a>5f$. Alors
l'ouvert ${\cal PW}$ des $m$ de ${\cal W}$ tels que $Z_m$ soit de
codimension exactement $f$ et le sous-espace associ\'{e} $Z_m^{*}$ de la
codimension attendue, \`{a} savoir $4f$, est non vide, irr\'{e}ductible, de
la codimension attendue $f(4b-10a+f)$ dans ${\cal M}$. 

De plus, la matrice
g\'{e}n\'{e}rique $\mu $ de ${\cal PW}$ repr\'{e}sente le morphisme
g\'{e}n\'{e}rique de $B$ vers le sous-espace $Z^{*}$ (de codimension $4f$)
de $A\otimes V$ associ\'{e} au sous-espace g\'{e}n\'{e}rique $Z$ de
codimension $f$ de $A\otimes S^2V$. En outre l'image de $\mu (1)$ est $Z$.
\end{proposition}

\noindent 
{\em D\'{e}monstration :} Soit $PG$ l'ouvert de la grassmannienne des
sous-espaces de codimension $f$ de $A\otimes S^2V$ qui sont de $V^{*}$-rang
maximum (lemme 3.2), et soit $Q$ l'ensemble des couples $(m,Z)$ de 
${\cal M}\times PG$ tels que $Z$ contienne l'image de $m(1)$. 
En projetant sur $PG$ et en
utilisant le transport d'\'{e}quations (proposition 2.1), 
on voit que c'est une
vari\'{e}t\'{e} lisse connexe non vide de dimension 
$d \=(10a-f)f+b(4a-4f)=4ab-f(4b-10a+f)$, c'est-\`a-dire de la dimension
attendue pour $W$. 

Soit $PQ$ l'ouvert de $Q$ o\`{u} $Z$ est
l'image de $m(1)$. Alors la premi\`{e}re projection induit un isomorphisme
de $PQ$ sur ${\cal PW}$. Il ne nous reste donc qu'\`{a} montrer
que $PQ$ est non vide, en  majorant par $d-1$ la
dimension du sch\'{e}ma $R$ des triplets $(m,Z,Z^{\prime })$ avec $(m,Z)$
dans $Q$, $Z^{\prime }$ hyperplan de $Z$ et l'image de $m(1)$ contenue dans 
$Z^{\prime }$. Pour cela, on stratifie $R$ par le $Z$-rang $r$ 
de $Z^{\prime }$, et l'on calcule la dimension des strates $R_r$ 
correspondantes par projection sur les deux derniers
facteurs. Gr\^{a}ce \`{a} la proposition 3.4, 
on obtient les majorations suivantes  :
$$\begin{array}{rcl}
\dim R_4 & \leq  & (10a-f)f+10a-f-1+b(4a-4f-4), \\
\dim R_3 & \leq  & (10a-f)f+ 6a+3f+2+b(4a-4f-3) \\
\dim R_2 & \leq  & (10a-f)f+\max(3a+3,-a+6f+3)+b(4a-4f-2) \\
\dim R_1 & \leq  & (10a-f)f+\max(-2a+11f+2,-a+4f+4,a+2)+b(4a-4f-1).
\end{array}$$
Ces nombres sont major\'{e}s par $d-1$ d\`{e}s que 
$4b\geq 10a-f,\; 3b\geq 6a+3f+3, \; 2b\geq 3a+4,\;
2b\geq -a+6f+4, \; b\geq -2a+11f+3, \; b\geq -a+4f+5$ et
$b\geq a+3$. Et l'on v\'erifie facilement que toutes ces 
in\'egalit\'es d\'ecoulent de nos hypoth\`eses. \qed

\medskip On arrive au r\'esultat principal de cette section, annonc\'e 
dans l'introduction (Th\'eor\`eme 1.1), et qui \'enonce l'existence
de la composante pr\'edominante du lieu de Brill-Noether $W$. On rappelle
que $W_P$ d\'esigne
l'ouvert de $W$ des
matrices $m$ telles que l'image de 
$$m(1)  : B\otimes V\longrightarrow A\otimes S^2V$$
soit de codimension exactement $f$, et o\`u le
morphisme $m \rightarrow Im(m(1))$, \`a valeurs dans la grassmannienne
des sous-espaces de codimension $f$ de $A \otimes S^2V$, est dominant.

\begin{theorem}
Soient $a$ et $b$ deux entiers positifs avec
$
5a/2 \leq b \leq 4a$, 
et soit $f$ positif. Pour $f$ sup\'erieur \`a $a - b/4$, $W_P$ est
contenu dans $W^0[0]$. Dans le cas
contraire ($f \leq a - b/4$), $W_P$ est contenu dans une unique
composante irr\'eductible $PW$ de $W$, qui a la codimension attendue
$f(4b - 10 a + f)$, et
dont le point
g\'en\'erique $\mu$ a les propri\'et\'es suivantes:
\begin{enumerate}
\item $\mu $ repr\'esente le morphisme g\'{e}n\'{e}rique de $B$
vers le sous-espace $Z^{*}$ (de codimension $4f$) de $A\otimes V$
associ\'{e} au sous-espace g\'{e}n\'{e}rique $Z$ de codimension $f$ de $%
A\otimes S^2V$; l'image de $\mu (1)$ est $Z$.
\item pour $d>1,\mu (d) :B\otimes S^dV\rightarrow A\otimes S^{d+1}V$ est
surjectif; donc \underline{$\mu $} est surjectif et son
noyau $E_\mu$ est localement libre de rang $b-a$.
\item
Le fibr\'{e} $E_\mu$ a la cohomologie attendue, \`{a} savoir  :
\begin{enumerate}
\item pour $k\neq 1$ l'un au plus des groupes $H^i(E_\mu (k))$ est non nul;
\item $h^1(E_\mu (1))=f$, ce qui d\'{e}termine tous les $h^i(E_\mu (1))$.
\end{enumerate}\end{enumerate}\end{theorem}

\noindent {\em D\'{e}monstration :} 
Supposons d'abord $a - b/4 < f$ et soit $\mu$ 
un point g\'en\'erique de $W_P$. L'image de $\mu(1)$
est contenue dans un sous-espace g\'en\'erique de codimension $f$ de $A \otimes 
S^2V$, donc aussi dans un sous-espace 
g\'en\'erique de codimension $f'$ de $A \otimes 
S^2V$, avec $a - b/4 < f' \leq a$. D'apr\`es le lemme 3.2, l'image de 
$\mu$ est alors contenue dans un sous-espace de codimension $4f'$ de
$A \otimes V$. De ce fait $\mu$ n'est pas injectif et appartient donc \`a
$W^0[0]$.

Supposons maintenant $f \leq a - b/4$, ce qui implique
$5f < 2a$. La proposition pr\'ec\'edente 
 fournit un ouvert irr\'{e}ductible ${\cal PW}$ de ${\cal W}$
 qui contient un ouvert dense de $W_P$ et
 caract\'{e}ris\'{e} par le point $1$. Nous allons v\'erifier
que son point g\'{e}n\'{e}rique $\mu $ v\'erifie les points $2$ et $3$.
Cela impliquera que la trace de ${\cal PW}$ sur $M$
est non vide et d\'{e}finit la composante $PW$ cherch\'ee. 

Montrons donc le point $2$. Comme l'image de $\mu (1)$ est
$Z$, celle de $\mu (d)$ est aussi l'image de $Z\otimes S^{d-1}V$. 
Il nous suffit donc de prouver que si $Y$ est le sous-espace 
g\'{e}n\'{e}rique de $A\otimes S^2V$ de dimension $9a<10a-f$, 
le morphisme naturel $Y\otimes S^{d-1}V\rightarrow
A\otimes S^{d+1}V$ est surjectif. On se ram\`{e}ne ainsi 
au cas $a=1$, que l'on traite directement en prenant $Y$ engendr\'{e} 
par neuf mon\^{o}mes dont les carr\'{e}s.

Ainsi $\mu $ d\'{e}termine un morphisme de modules gradu\'{e}s dont le coker
a longueur finie, par cons\'{e}quent \underline{$\mu $} est
surjectif, et son noyau $E_{\mu}$ est bien un fibr\'e de rang $b-a$.

On calcule alors, pour \'etablir le point $3$, 
les groupes de cohomologie de $E_{\mu}$ \`{a} l'aide de la suite
exacte
$$0\longrightarrow E_\mu \longrightarrow B\otimes {\cal O}
\longrightarrow A\otimes{\cal O}(1)\longrightarrow 0.$$ 
L'\'{e}nonc\'{e} (a) d\'{e}coule alors du point $2$ (pour $k \geq 2$), 
et du fait que $\mu$ 
est le morphisme g\'{e}n\'{e}rique de $B$ vers $Z^{*}$ (pour $k = 0$).
L'\'enonc\'e (b) r\'esulte
du fait que $Im(\mu (1))=Z$. \qed

\medskip\noindent {\bf Remarque 4.3}
La cohomologie de $E_\mu $ est donn\'{e}e par
le tableau suivant, o\`u $\chi (k)$ d\'esigne la caract\'eristique 
d'Euler de $E_\mu(k)$  :
\medskip\begin{center}
\begin{tabular}{||c|c|c|c|c|c|c|c||} \hline
$k$ & $\cdots$ & -2 & -1 & 0 & 1 & 2 & $\cdots$ \\ \hline
$h^0(E_\mu(k))$ & 0 & 0 & 0 & 0 & $4b-10a+f$ & $\chi (2)$ & $\cdots$ \\
$h^1(E_\mu(k))$ & 0 & 0 & $a$ & $4a-b$ & $f$ & 0 & 0 \\
$h^2(E_\mu(k))$ & 0 & 0 & 0 & 0 & 0 & 0 & 0 \\
$h^3(E_\mu(k))$ & $\cdots$ & $\chi (-2)$ & 0 & 0 & 0 & 0 & 0 \\ \hline
\end{tabular}\end{center}\medskip
Observons que dans le cas que nous consid\'{e}rons au \S 7,
avec $3a\geq b$ et $f=9a-3b+1$, on a $h^0(E_\mu (1))
\leq b-a+1={\rm rang}(E_\mu)+1.$
Ceci implique que $E_\mu (1)$ {\em n'est pas engendr\'{e} par ses 
sections globales}  :  si c'\'etait le cas, on aurait en effet une suite 
exacte 
$$0\longrightarrow L \longrightarrow S
\longrightarrow E_{\mu}(1)\longrightarrow 0,$$ 
avec $S$ fibr\'e trivial et $L$ fibr\'e de rang un, ce qui impliquerait 
$h^1(E_{\mu})=h^1(E_{\mu}(1))=0$. 

\section{Stratification par le $(Z,H)$-rang}

Soit $Z$ le sous-espace g\'en\'erique de codimension $f$ de $A\otimes S^2V$
et $H$ un hyperplan de $V$ (d\'efini sur le corps 
de d\'efinition de $Z$) .

\begin{lemma}
Si $3a\ge f+3$, $Z$ est transverse \`{a} $A\otimes H.V$.
\end{lemma}

\noindent {\em D\'emonstration} :
Les sous-espaces de codimension $f$ de $A\otimes S^2V$ non
transverses \`{a} $A\otimes H'.V$, o\`u $H'$ est un hyperplan fix\'e de
$V$,  forment un cycle de Schubert de codimension
$3a-f+1$ dans la grassmannienne correspondante. Ceux qui sont
non-transverses \`{a} l'un au moins des $A\otimes H'.V$ forment 
donc un cycle de codimension au moins $3a-f-2$, 
et ce nombre est strictement positif.\qed

\medskip\noindent
On note d\'esormais $Z'$ l'intersection (transverse, comme on vient de
le voir) de $Z$ et $A\otimes H.V$. On introduit un nouveau rang
auxiliaire, le
$(Z,H)$-rang, \`a l'aide duquel on va stratifier successivement
l'espace projectif des hyperplans de $Z'$ et la grassmannienne de
ses sous-espaces de codimension deux.
 
\begin{definition} (du $(Z,H)$-rang)
Soit $T$ un sous-espace de $Z^{\prime } \=Z\cap A\otimes H.V$. 
On pose $F^{\prime } \=(A\otimes H.V)/T,\; F \=(A\otimes S^2V)/Z$, 
et l'on note 
$$f_{Z,T}^{*} :A\otimes V\longrightarrow 
(F\otimes V^{\vee})\oplus (F^{\prime }\otimes H^{\vee})$$ 
le morphisme correspondant (cf. 2.3). On d\'{e}finit alors le 
$(Z,H)$-rang de $T$ comme le rang de cette application $f_{Z,T}^{*}$, 
diminu\'{e} du $V^{*}$-rang de $Z$. \end{definition}

On notera $G_r$ le sous-sch\'ema des hyperplans de $Z'$ de
$(Z,H)$-rang au plus $r$,
sous-sch\'ema dont on va majorer la dimension.

\begin{proposition} 
Supposons  $3a>11f$. Alors $G_0$ est vide, et l'on a
$$\begin{array}{lcl}
{\rm codim}\; G_1 & \geq & \min (8a-8f-2,7a-f-2), \\
{\rm codim}\; G_2 & \geq & 4a-4f-2.
\end{array}$$ \end{proposition} 

\noindent
{\em D\'{e}monstration :} Soit $\Phi$ un morphisme surjectif de $A\otimes
S^2V$ vers $k^f$ dont $Z$ est le noyau, et $t_H$ sa restriction \`{a} 
$A\otimes H.V$. Au lieu de $G_r$ il nous suffit de stratifier le sch\'{e}ma,
de dimension $9a$, des (noyaux des) morphismes surjectifs $g$ de $A\otimes
H.V$ vers $k^{f+1}$ de la forme $(t_H,-)$. Choisissons  des bases 
$(e_j) $ sur $A$ et $(v_p)$ sur $V$ avec $H=\left\{ v_4=0\right\} $. La
matrice sym\'{e}trique $X$ de $t$ s'\'{e}crit $x_{pqjs}$ avec $1\leq p,q\leq
4,1\leq j\leq a $ et $1\leq s\leq f$. Celle, partiellement sym\'{e}trique,
de $g$ s'\'{e}crit $g_{pqjs}$ avec $1\leq p\leq 3,1\leq q\leq 4,1\leq j\leq
a $ et $1\leq s\leq f+1$, et $g_{pqjs}=x_{pqjs}$ pour $s\leq f$. 
On pose $g_{pqj}=g_{pqj,f+1} $. Le $(Z,H)$-rang de $\ker (g)$ est donc le
surcro\^{i}t de rang qu'apportent les ''lignes'' $g_{p..}$ ($p\leq 3$) aux
''lignes'' $x_{p..s}$ ($p\leq 4,s\leq f$) de $t$. Comme plus haut, les
majorations de dimension se font localement.

\medskip\noindent {\bf Cas $r=0$}.
On proc\`ede comme pour la proposition 3.4, 
\`a ceci pr\`es que  $G$ et $C$ sont ici
des matrices $3\times 4$. On \'ecrit $G= (G_0,g)$, $C= (C_0,c)$
et $$X= \left ( \begin{array}{cc} X_0 & x \\ x^t & \xi
\end{array} \right ).$$
On calcule donc les dimensions des strates et des fibres du
sch\'ema $\Gamma$ des couples $(C, X)$ avec $CX$ partiellement
sym\'etrique, et d'abord dans le cas
$f= 1$. 
Si $P\in Gl(3,\CC)$, soit $i(P)$ son image dans $Gl(4,\CC)$. Alors 
$G= CX$ est partiellement sym\'etrique si et seulement si $PGi(P)^t= 
(PCi(P)^{-1})(i(P)Xi(P)^t)$ l'est aussi. Comme $PCi(P)^{-1}
= (PC_0P^{-1},P(c))$, on peut supposer que $C_0$ est
sous forme de Jordan. La matrice $G$ est 
partiellement sym\'etrique si et seulement si la matrice 
$C_0X_0+c.x^t$ est sym\'etrique, d'o\`u un syst\`eme
de trois \'equations en les coefficients de $X$. On stratifie donc
l'espace des matrices $C$ par la forme de Jordan de $C_0$ et la
forme de $c$ de fa\c{c}on \`a pouvoir calculer le rang $r$ du
syst\`eme pr\'ec\'edent. On obtient le tableau suivant, dans lequel
$O$ d\'esigne la dimension des strates et $S=10-r$ est la dimension de
l'espace des solutions correspondantes :
\medskip\begin{center}
\begin{tabular}{||c|c|c|c|c||} \hline
$J$ & $c$ & $r$ & $O$ & $S$ \\ \hline
$1|1|1$ &  & 3 & 12 & 7 \\
$11|1$ & $(0,0,*)$ & 2 & 7 & 8 \\
     & {\rm autre} & 3 & $<12$ & 7 \\
$111$ & $\neq$ 0 & 2 & 7 & 8\\ 
    &  0   & 0 & $g$ {\rm non surjectif} & \\
$2|1$ & & 3 & $<12$ & 7 \\
$21$ & $(*,0,0)$ & 2 & 6 & 8 \\
   & {\rm autre} & 3 & $<12$ & 7 \\
3 &  & 3 & $<12$ & 7 \\ \hline
\end{tabular}\end{center}\medskip

Lorsque $f$ est quelconque, on consid\`ere $C$ comme une suite
$C= (C^1,\ldots ,C^f)$ de $f$ matrices $3\times 4$, et l'on 
stratifie selon le maximum  $O(C)$ des $O(C^i)$, tandis
qu'on note $S(C)$ le minimum correspondant des $S(C^i)$.
Les dimensions $O_f(C)$ et $S_f(C)$ des strates 
et des espaces de solutions correspondants v\'erifient \'evidemment
$$O_f=  f\times O,\quad S_f\leq S+10(f-1).$$
On observe maintenant que les fibres de la restriction \`a $\Gamma$ de
la premi\`ere projection sont de dimension au moins $f$,
\`a cause de l'invariance par
$c_{pqs}\rightarrow c_{pqs}+\gamma_s\delta _{pq}$.
On constate alors que toute les strates v\'erifient
$(10-S)a>(O-1)f$, ce dont on d\'eduit l'in\'egalit\'e $10af>aS_f+O_f-f$,
qui signifie qu'aucune strate ne domine l'espace des matrices $X$.

\medskip\noindent {\bf Cas $r=1$}. 
Quitte \`{a} changer de base dans $H$, on peut supposer que le
surcro\^{\i}t de rang est d\^u \`a $g_{1..}$,
d'o\`u des relations de d\'ependance
$$\begin{array}{rcl}
g_{2..} & =  & \alpha _2g_{1..}+\sum_{rs}c_{2rs}x_{r..s}, \\
g_{3..} & =  & \alpha _3g_{1..}+\sum_{rs}c_{3rs}x_{r..s}.
\end{array}$$
On doit tenir compte de l'unique relation de sym\'etrie $g_{23.}=g_{32.}$,
c'est-\`a-dire 
$$\sum_{rs}c_{3rs}(x_{r2.s}-\alpha _2x_{r1.s})
= \sum_{rs}c_{2rs}(x_{r3.s}-\alpha _3x_{r1.s}).$$
L\`a o\`u cette relation est non-triviale, 
la strate correspondante de $\Gamma_1$ est de dimension
$2$ (pour les $\alpha$) plus $8f$ (pour $c$) plus $2a$ (pour $g_{11.}
et g_{14.}$) plus $a(9f-1)$ (pour les $x$ autres que $x_{44}$), et 
donc sa fibre \`a $x$ fix\'e  est de dimension au plus
$a+8f+2$. La codimension de cette fibre est donc au moins $8a-8f-2$.
Si la relation est triviale, c'est que
$c_{pqs}= 0$ pour $2\leq p\neq q\leq 4$ et $p\leq 3$, 
que $c_{22s}= c_{33s} \= c_s$ et que $c_{p1s}= -\alpha _pc_s$
pour $p= 2,3$. La strate correspondante de $\Gamma_1$ est de dimension
$2$ (pour les $\alpha$) plus $f$ (pour les $c_s$) plus $2a$ (pour $g_{11.}$
et $g_{14.}$) plus $9af$ (pour les $x$ autres que $x_{44}$), et 
donc sa fibre \`a $x$ fix\'e est de dimension au plus
$2a+f+2$. La codimension de cette fibre est donc au moins $7a-f-2$.

\medskip\noindent {\bf Cas $r=2$}.
On peut de m\^{e}me supposer que le surcro\^{\i}t de rang est
d\^{u} \`{a} $g_{1..}$ et $g_{2..}$. L'ouvert correspondant est
param\^{e}tr\'{e} par $g_{11.},g_{12.},g_{14.},g_{22.},g_{24.}$ et les
coefficients exprimant $g_{3..}$ en fonction lin\'{e}aire des 
$x_{p..s}$ et de $g_{1..}$ et de $g_{2..}$. On majore ainsi la dimension 
par $5a+4f+2$.\qed

\medskip\noindent {\bf Remarque 5.4} Si $3a<11f$, alors $G_0$ est 
non vide. En effet, reprenons les notations de la d\'emonstration 
pr\'ece\'edente, et raisonnons dans l'espace des matrices $g$ non
n\'ecessairement sym\'etriques, espace qui est de dimension $12a$.
Parmi ces matrices, celles qui ne d\'efinissent pas un hyperplan 
forment un sous-espace de dimension $f$. Celles pour lesquelles 
$r=0$, c'est-\`a-dire pour lesquelles on a une relation de 
d\'ependance de la forme 
$$g_{p..}=\sum_{qs}c_{pqs}x_{q..s},$$
forment un sous-espace de dimension $12f$. Parmi celles-ci,
celles qui sont sym\'etriques forment un sous-espace de 
dimension au moins $12f-3a$, puisque le sous-espace des matrices 
sym\'etriques est de codimension $3a$. Si $12f-3a>f$, certaines 
de ces matrices de\'finissent n\'ecessairement des hyperplans de 
$Z$ qui sont des points de $G_0$. 
 
\bigskip Notons maintenant $G'_r$ le sous-sch\'ema des sous-espaces de 
codimension deux de $Z'$ de $(Z,H)$-rang au plus $r$, et majorons sa dimension.

\addtocounter{lemma}{1}
\begin{proposition} 
Supposons que $2a>5f$. Alors $G'_0$ est vide, et on a
$$\begin{array}{lcl}
{\rm codim}\; G'_1 & \geq & \min (20a-20f-5,19a-12f-7,18a-7f-6,17a-2f-5), \\
{\rm codim}\; G'_2 & \geq & \min (16a-16f-8,15a-9f-7,14a-2f-4), \\
{\rm codim}\; G'_3 & \geq & \min (12a-12f-8,11a-5f-5), \\
{\rm codim}\; G'_4 & \geq & \min (8a-8f-8,7a-f-3), \\
{\rm codim}\; G'_5 & \geq & 4a-4f-5.
\end{array}$$
\end{proposition} 

\noindent {\em D\'{e}monstration :}
On stratifie cette fois le sch\'{e}ma, de dimension $18a$, 
des morphismes $g$ de $A \otimes H.V$ vers 
$k^{f+2}$ de la forme $(f_H,-)$. On garde les m\^emes notations 
en posant de plus $g'_{pqj }= g_{pqj,f+2 }$.
On voit ici $X$ comme un morphisme de $V\otimes F^{\vee}$ vers
$A^{\vee}\otimes V^{\vee}$,
$g$ et $g'$ comme deux morphismes de $H$ vers le m\^eme but,
$c$ et $c'$ comme deux morphismes de $H$ vers
$V\otimes F^{\vee}$, 
et $h$, $k$, ... comme des vecteurs de $A^{\vee}\otimes V^{\vee}$.

\medskip\noindent {\bf Cas $r=0$}. Cela d\'ecoule de la proposition 
pr\'ec\'edente.

\medskip\noindent {\bf Cas $r=1$}. On \'ecrit, pour $1\leq p\leq 3$, 
$$\begin{array}{rcl}
g_{p..} & =  & \alpha _ph_{..}+\sum_{rs}c_{prs}x_{r..s}, \\
g'_{p..} & =  & \alpha' _ph_{..}+\sum_{rs}c'_{prs}x_{r..s},
\end{array}$$
et on estime la dimension du sch\'ema des 
$(x,h,\alpha ,\alpha ',c,c')$ pour lesquels la contrainte
de sym\'etrie est respect\'ee, en projetant sur 
$(\alpha ,\alpha ',c,c')$. Notons que les fibres de la projection 
de $(x,h,\alpha ,\alpha ',c,c')$ sur $(x,g,g')$ sont de dimension 
$4l+1$, \`a cause de l'invariance par la transformation 
$$\begin{array}{lcl}
\alpha_p & \mapsto & m^{-1}\alpha_p, \\ 
\alpha'_p & \mapsto & m^{-1}\alpha'_p, \\ 
h_{..} & \mapsto & mh_{..}+\sum_{rs}n_{rs}x_{r..s}, \\
c_{prs} & \mapsto & c_{prs}-\alpha_pn_{rs}, \\
c'_{prs} & \mapsto & c'_{prs}-\alpha'_pn_{rs}.
\end{array}$$ 
Les six relations de sym\'etrie 
s'\'ecrivent, $A$ \'etant mis en facteur
$$\begin{array}{rcl}
\alpha _ph_q-\alpha _qh_p & =  & \sum_{rs}(c_{qrs}x_{prs}
-c_{prs}x_{qrs}), \\
\alpha' _ph_q-\alpha' _qh_p & =  & \sum_{rs}(c'_{qrs}x_{prs}
-c'_{prs}x_{qrs}),
\end{array}$$
avec $1\leq p<q\leq 3$.
Ce syst\`eme s'\'ecrit sous la forme 
$$h\wedge \alpha = C(x),\quad {\rm et}\quad h\wedge \alpha '=C'(x),$$
o\`u par exemple $h\wedge \alpha$ 
d\'esigne l'image naturelle  de $h\otimes \alpha$ dans
$A^{\vee}\otimes \Lambda^2 H^{\vee}$; tandis que $C(x)$ est l'image naturelle 
dans $A^{\vee}\otimes \Lambda^2 H^{\vee}$ du compos\'e $X\circ c$.
On remarque alors que $\alpha $ et $\alpha '$ ne peuvent
\^etre parall\`eles, \`a cause de la premi\`ere assertion 
de la proposition 3.4. Le syst\`eme pr\'ec\'edent
d\'etermine donc la projection de $h$ 
dans $A^{\vee}\otimes H^{\vee}$ sous r\'eserve que soient v\'erifi\'ees les
trois conditions de compatibilit\'e
$$C(x)\wedge \alpha = C'(x)\wedge \alpha' = C(x)\wedge \alpha' +
C'(x)\wedge \alpha = 0.$$
On stratifie donc par le rang de ce syst\`eme en $x$, en choisissant
une base commen\c{c}ant par $\alpha$ et $\alpha'$. Le syst\`eme devient
$$\begin{array}{rcl}
\sum_{rs}(c_{2rs}x_{3rs}-c_{3rs}x_{2rs}) & =  & 0, \\
\sum_{rs}(c'_{1rs}x_{3rs}-c'_{3rs}x_{1rs}) & =  & 0, \\
\sum_{rs}(c_{1rs}x_{3rs}-c_{3rs}x_{1rs}) & =  & 
\sum_{rs}(c'_{2rs}x_{3rs}-c'_{3rs}x_{2rs}).
\end{array}$$
Si ce syst\`eme est de rang z\'ero, on obtient $c_{pqs}= 0$ pour 
$p= 2,3$ et $q\neq p$, $c'_{pqs}= 0$ pour $p= 1,3$ et $q\neq p$, ainsi
que $c_{22}= c_{33}$,  $c'_{11}= c'_{33}$, 
$c_{13s}= c'_{23s}$, $c_{14s}= c'_{24s}$, $c_{12s}= c'_{22s}-c'_{33s}$
et $c'_{21s}= c_{11s}-c_{33s}$. 
Cette strate est de dimension $6$ (pour $\alpha$ et $\alpha'$) plus
$6f$ (pour $c_{11}, c_{22}, c'_{11}, c'_{22}, c'_{23} et c'_{24}$),
plus $10af$ (pour les $x$) plus $a$ (pour les composantes de $h$
qui ne sont pas d\'etermin\'ees) moins $4f+1$ (pour les fibres de la
projection sur $(g,g')$), soit de codimension $17a-2f-5$.

S'il est de rang un sans qu'aucune \'equation ne soit identiquement 
nulle, elles doivent \^etre proportionnelles \`a une combinaison
de $x_{12},x_{13},x_{23},x_{33}$ et $x_{34}$, en fait des quatre
derniers seulement, et 
$c$ et $c'$ se calculent en fonction des $c_{1.}$, des $c_{2.}$ et de
$c'_{11}$. La strate correspondante est donc de dimension $6$ (pour
$\alpha$ et $\alpha'$) plus $2$ (pour les d\'ependances entre les
trois \'equations) plus $9f$ (pour $c$ et $c'$),
plus $a(10f-1)$ (pour les $x$) plus $a$ (pour les composantes de $h$
qui ne sont pas d\'etermin\'ees) moins $4f+1$ (pour les fibres de la
projection sur $(g,g')$), soit de codimension $18a -5f-7$.
Si une des trois \'equations est identiquement nulle, on v\'erifie que
les strates correspondantes sont  de dimension au plus $6$ (pour
$\alpha$ et $\alpha'$) plus $1$ (pour la d\'ependance entre les
deux \'equations restantes) plus
$10f$ (pour $c$ et $c'$, l\`a est la v\'erification),
plus $a(10f-1)$ (pour les $x$) plus $a$ (pour les composantes de $h$
qui ne sont pas d\'etermin\'ees) moins $1$ (pour les fibres de la
projection sur $(g,g')$); soit de codimension au moins $18a-10f-6$.

Passons au cas o\`u le syst\`eme est de rang deux. Si $c_{34.}$ ou 
$c'_{34.}$ est non nul, les coordonn\'ees $(a :b :c)$ de la relation
v\'erifient $ab+c^2= 0$, et cette strate est de dimension $6$ 
(pour $\alpha$ et $\alpha'$) plus $1$ (pour la relation pr\'ec\'edente), 
plus $17f$ (pour $c$ et $c'$; par exemple si $c'_{34.}$ est non nul,
$c_{32}, c_{34}, c_{31}, c_{22}, c_{23}, c_{21}$ et $c_{24}$ se
calculent en fonction des autres), plus $a(10f-2)$ (pour les $x$) 
plus $a$ (pour les composantes de $h$ qui ne sont pas d\'etermin\'ees)
moins $4f+1$ (pour les fibres de la projection  sur $g,g'$), 
soit de codimension $19a-13f-6$.
Si maintenant $c_{34.}$ et $c'_{34.}$ sont nuls, on distingue suivant
l'\'equation dont la relation de d\'ependance permet de calculer les 
coefficients, et on constate dans chaque cas qu'on peut calculer six
nouveaux coefficients en fonction des autres. Les strates
correspondantes sont alors de codimension au moins $19a-8f-7$.

Enfin la strate o\`u le rang est trois est de dimension $6$ 
(pour $\alpha$ et $\alpha'$), plus $24f$ (pour $c$ et $c'$),
plus $a(10f-3)$ (pour les $x$) plus $a$ (pour les composantes de $h$
qui ne sont pas d\'etermin\'ees) moins $4f+1$ (pour les fibres de la
projection sur $(g,g')$), soit de codimension $20a-20f-5$.

\medskip\noindent {\bf Cas $r=2$}.
Comme ci-dessus, on \'ecrit
$$\begin{array}{rcl}
g_{p..} & =  & \alpha _ph_{..}+\beta _pk_{..}+\sum_{rs}c_{prs}x_{r..s}, \\
g'_{p..} & =  & \alpha'_ph_{..}+\beta'_pk_{..}+\sum_{rs}c'_{prs}x_{r..s}, 
\end{array}$$
avec $1\leq p\leq 3$. Ici, les fibres de la projection 
de $(x,h,k,\alpha ,\alpha ',\beta, \beta',c,c')$ sur $(x,g,g')$ 
sont de dimension $8f+4$, \`a cause de l'invariance par 
l'action de $Gl_2$ sur $h,k,\alpha ,\alpha ',\beta, \beta'$, sur
laquelle on va revenir, et par la transformation 
$$\begin{array}{lcl}
h_{..} & \mapsto & h_{..}+\sum_{rs}l_{rs}x_{r..s}, \\
k_{..} & \mapsto & k_{..}+\sum_{rs}n_{rs}x_{r..s}, \\
c_{prs} & \mapsto & c_{prs}-\alpha_pl_{rs}, \\
c'_{prs} & \mapsto & c'_{prs}-\alpha'_pn_{rs}.
\end{array}$$ 
Les relations de sym\'etrie donnent, $C$ et $C'$ \'etant comme ci-dessus,
$$\begin{array}{rcl}
\alpha\wedge h +\beta\wedge k  & =  & C(x), \\
\alpha'\wedge h +\beta'\wedge k  & =  & C'(x).
\end{array}$$
Et on estime la dimension du sch\'ema des 
$(x, h, k,\alpha ,\alpha ', \beta, \beta', c, c')$ pour lesquels la contrainte
de sym\'etrie est respect\'ee, en projetant sur le sch\'ema des
$(h, k,\alpha ,\alpha ', \beta, \beta', c, c')$.

\medskip\noindent
{\sc Premi\`ere \'etape  : r\'eduction de $\alpha ,\alpha ', \beta, \beta'$.}

Observons tout d'abord que $g$ et $g'$ sont inchang\'es si l'on compose
la matrice $2\times 2$ de trivecteurs $(\alpha ,\beta |
\alpha',\beta')$ \`a droite par une
matrice $P\in Gl_2$, tout en appliquant $P^{-1}$ au vecteur $(h,k)$.
On peut \'egalement changer de base de fa\c con \`a remplacer $g,g'$
par des combinaisons lin\'eaires, ce qui a pour effet de multiplier
la matrice $(\alpha ,\beta | \alpha' ,\beta')$ \`a gauche par une
nouvelle matrice $Q\in Gl_2$. 

\smallskip
Ceci permet de r\'eduire la matrice $(\alpha ,\beta | \alpha' ,\beta')$,
que l'on consid\`ere, dans une base donn\'ee, comme une famille de trois
matrices $A_1,A_2,A_3$. A cause de la premi\`ere assertion de la 
proposition 5.3, ces matrices ne sont pas
toutes nulles. On peut m\^eme supposer qu'elles ne sont pas toutes 
singuli\`eres, puisque si c'\'etait le cas dans toute base, 
elles seraient toutes les trois multiples d'une m\^eme matrice de rang
un, que l'on pourrait supposer \^etre $(1,0 | 0,0)$. Autrement dit, on
serait ramen\'e \`a $\beta= \alpha'= \beta'= 0$, et les $g'_{p..}$ seraient
combinaisons lin\'eaires des $x_{p..s}$, ce qu'\`a nouveau 
la proposition 5.3 exclut.

\smallskip
On peut donc supposer $A_1$ inversible, et m\^eme $A_1= Id$ gr\^ace \`a
l'action de $P$ et $Q$. On peut encore composer par ces matrices 
si $P= Q^{-1}$, donc r\'eduire $A_2$ \`a sa forme de Jordan.

\smallskip\noindent
{\em Premier cas}  : $A_2$ est une homoth\'etie. On peut alors aussi 
r\'eduire $A_3$ \`a sa forme canonique, d'o\`u trois sous-cas  :
\begin{enumerate}
\item $A_3$ est encore une homoth\'etie, et l'on est ramen\'e \`a 
$\alpha= \beta'$ et $\alpha'= \beta= 0$. Les orbites correspondantes
sous $Gl(4,\CC)\times Gl(4,\CC)$ sont de dimension $4$, 
cette strate est donc de dimension $6$.
\item $A_3$ n'est pas diagonalisable, auquel cas on est ramen\'e 
\`a $\alpha= \beta'$ ind\'ependant de $\beta$,  et $\alpha'= 0$. 
Les orbites sont de dimension $6$, cette strate est donc de dimension $8$.
\item $A_3$ a deux valeurs propres distinctes, et l'on est ramen\'e \`a
$\alpha'= \beta= 0$. Cette strate est de dimension $9$.
\end{enumerate}

\smallskip\noindent
{\em Deuxi\`eme cas}  : $A_2$ n'est pas diagonalisable. Si on lui
donne sa forme canonique $(m,1|0,m)$, on peut encore conjuguer 
$A_3= (a,b|c,d)$ 
par $P$ de la forme $(1,p|0,1)$, ce qui permet de supposer que 
$a= 0$ ou $c= 0$. D'o\`u deux sous-cas  :
\begin{enumerate}
\item si $c= 0$, on est ramen\'e \`a $\alpha'= 0$. On peut supposer
$a\neq d$, sans quoi on retrouve le premier cas. Les vecteurs
$\alpha,\beta$ et $\beta'$ sont alors ind\'ependants.
Les orbites sont de dimension $7$, cette strate est donc de dimension $11$.
\item si $c\neq 0$, $a= 0$, il existe une base  $\gamma,\gamma',\gamma''$ 
telle que $\alpha= \gamma+e\gamma'$, $\beta= \gamma'+b\gamma''$,
$\alpha'= \gamma''$ et $\beta'= \gamma+e\gamma'+d\gamma''$.
Cette strate est aussi de dimension $11$.
\end{enumerate}

\smallskip\noindent
{\em Troisi\`eme cas}  : $A_2$ a deux valeurs propres distinctes. 
Quitte \`a changer de base dans $H$, on peut supposer $A_3$ de la forme 
$(0,b|c,0)$. Si $c= 0$, on retrouve le deuxi\`eme cas, \`a moins 
d'avoir aussi $b= 0$, qui redonne le premier cas. On peut donc
supposer $c\neq 0$. Autrement dit, on peut supposer que
$\alpha'$ et $\beta$ sont collin\'eaires, tous deux non nuls, et 
que les vecteurs $\alpha,\beta,\beta'$ sont ind\'ependants.
Cette strate est la g\'en\'erique, de dimension $12$.

\medskip\noindent
{\sc Deuxi\`eme \'etape  : estimation des dimensions.}

\smallskip\noindent
{\em Premier cas}  : On reprend s\'epar\'ement les trois sous-cas.
\begin{enumerate}
\item Les \'equations de sym\'etrie se r\'eduisent \`a 
$\alpha\wedge h =  C(x)$ et $\alpha\wedge k  =  C'(x)$.
Comme $\alpha$ n'est pas nul, deux composantes de chacun des 
vecteurs $h$ et $k$ sont d\'etermin\'ees, sous r\'eserve que soient
v\'erifi\'ees les relations de compatibilit\'e 
$$C(x)\wedge \alpha =  C'(x)\wedge \alpha= 0.$$
Dans l'espace des couples $(c,c')$, ce syst\`eme est g\'en\'eriquement de
rang deux, de rang un en codimension $7f-1$, et de rang z\'ero en 
codimension $14f$. 

Selon les cas, on majore donc la dimension par $6$ 
(pour $\alpha, \alpha',\beta,\beta'$) plus  $24f$ 
(respectivement $17f+1$ et $10f$, pour $c$ et $c'$) 
plus $10af-2a$ (respectivement $10af-a$ et $10af$, pour $x$) 
plus $4a$ (pour les composantes non d\'etermin\'ees de $h$ et $k$) 
moins $8f+4$ (pour les fibres de la projection sur $(g,g')$).
Cette strate est donc de codimension minor\'ee par le minimum de
$16a-16f-2$, $15a-9f-3$ et $14a-2f-2$.

\item Ici, les \'equations de sym\'etrie donnent
$\alpha\wedge h +\beta\wedge k=  C(x)$ et $\alpha\wedge k  =  C'(x)$.
Comme $\alpha$ et $\beta$ sont ind\'ependants, ce syst\`eme d\'etermine 
la projection de $k$ 
dans $A^{\vee}\otimes H^{\vee}$, sous r\'eserve des conditions de compatibilit\'e 
$$\beta\wedge \alpha\wedge h =  C(x)\wedge \beta, \quad
C'(x)\wedge \alpha =  C(x)\wedge \alpha +C'(x)\wedge \beta= 0.$$
La premi\`ere de ces \'equations d\'etermine une composante de $h$.
Les deux derni\`eres donnent un syst\`eme en $x$ qui, 
dans l'espace des couples $(c,c')$, est g\'en\'eriquement de
rang deux, de rang un en codimension $9f-1$, et de rang z\'ero en 
codimension $14f$. On minore donc la codimension de la strate 
correspondante en $(g,g')$, par le minimum de
$16a-16f-4$, $15a-7f-5$ et $14a-2f-4$.

\item De la m\^eme fa\c con qu'au sous-cas pr\'ec\'edent, 
on majore la codimension de la strate 
en $(g,g')$, par le minimum de $16a-16f-5$, $15a-7f-6$ et $14a-2f-5$.
\end{enumerate}

\smallskip\noindent
{\em Deuxi\`eme cas}  : On a deux sous-cas.
\begin{enumerate}
\item Les \'equations de sym\'etrie se r\'eduisent \`a 
$$\alpha\wedge h +\beta\wedge k=  C(x), \quad \beta'\wedge k  =  C'(x).$$
Comme $\beta$ et $\beta'$ sont ind\'ependants, elles  d\'eterminent 
la projection de $k$ dans $A^{\vee}\otimes H^{\vee}$,
sous r\'eserve des conditions de compatibilit\'e
$$\begin{array}{rcl}
\beta\wedge\alpha\wedge h & =  & C(x)\wedge\beta, \\
\beta'\wedge\alpha\wedge h & =  & C(x)\wedge\beta'+C'(x)\wedge\beta, \\
C'(x)\wedge\beta'  & =  & 0.
\end{array}$$
Les deux premi\`eres \'equations d\'eterminent deux composantes de $h$,
puisque les vecteurs $\alpha,\beta$ et $\beta'$ sont ind\'ependants, 
donc $\beta\wedge\alpha$ et $\beta'\wedge\alpha$ aussi.
La derni\`ere \'equation est identiquement nulle en $x$, en codimension 
$7f$ dans l'espace des $c'$. On en d\'eduit que la strate correspondante
est de codimension minor\'ee par le minimum de $16a-16f-7$ et $15a-9f-7$.

\item  Ici, les conditions de sym\'etrie d\'eterminent la projection 
de $k$ dans $A^{\vee}\otimes H^{\vee}$, et les  relations de 
compatibilit\'e sont 
$$\begin{array}{rcl}
\beta\wedge\alpha\wedge h & =  & C(x)\wedge\beta, \\
\beta'\wedge\alpha'\wedge h & =  & C'(x)\wedge\beta', \\
(\beta'\wedge\alpha+\beta\wedge\alpha')\wedge h & =  & 
C(x)\wedge\beta'+C'(x)\wedge\beta.
\end{array}$$
Les trois vecteurs $\beta\wedge\alpha$, $\beta'\wedge\alpha'$
et $\beta\wedge\alpha'+\beta'\wedge\alpha$ sont ind\'ependants, donc 
la projection de $h$ dans $A^{\vee}\otimes H^{\vee}$ 
est aussi d\'etermin\'ee. On en d\'eduit que la strate correspondante
en $(g,g')$ est de codimension $16a-16f-7$.
\end{enumerate}

\smallskip\noindent
{\em Troisi\`eme cas}  : M\^eme estimation qu'au sous-cas
pr\'ec\'edent, si ce n'est que la dimension de la strate de
$(\alpha,\beta|\alpha',\beta')$ est sup\'erieure d'une unit\'e.
D'o\`u une codimension $16a-16f-8$ au moins.

\medskip\noindent {\bf Cas $r=3$.}
On doit a priori envisager deux cas selon la r\'{e}partition
$(3+0$ et $2+1)$ du surcro\^{\i}t de rang entre $g_{...}$ et 
$g_{...}^{^{\prime }}$. 
Montrons tout d'abord que le premier cas ne peut pas se
pr\'{e}senter. En effet, quitte \`{a} changer de
coordonn\'{e}es dans $k^{f+2}$, on pourrait supposer que 
$g_{...}^{^{\prime }}$
est combinaison lin\'{e}aire des $x_{p..s}$, ce qui n'est pas possible,
comme on l'a vu en 5.2, cas a).

Supposons donc que le surcro\^{\i}t de rang 
d\^{u} \`a $g$ est \'egal \`a deux (disons qu'il est d\^{u} \`a $g_{1..}$ et 
$g_{2..}$), et que celui d\^{u} en sus \`a $g'$ est \'egal \`a un (disons
qu'il est d\^{u} \`a $g'_{1..}$). On doit alors avoir des relations de 
d\'ependance de la forme
$$\begin{array}{rcl}
g_{3..} & =  & ag_{1..}+bg_{2..}+\sum_{ps}l_{ps}x_{p..s}, \\
g'_{2..} & =  & cg_{1..}+dg_{2..}+eg'_{1..}+\sum_{ps}m_{ps}x_{p..s}, \\
g'_{3..} & =  & fg_{1..}+hg_{2..}+kg'_{1..}+\sum_{ps}n_{ps}x_{p..s}.
\end{array}$$
L'ouvert correspondant est donc param\'etr\'e par $g_{11.}$, $g_{12.}$,
$g_{22.}$, $g_{14.}$, $g_{24.}$, $g'_{11.}$,$g'_{14.}$ et par 
les coefficients $l_{ps}, m_{ps}, n_{ps}$, li\'es par la seule relation 
de sym\'etrie $g'_{23.}= g'_{32.}$. Celle-ci se traduit par l'identit\'e
$$\begin{array}{rll}
Ag_{11.}+Bg_{12.}+Cg_{22.} & 
+\sum_{ps}m_{ps}x_{p3.s}+c\sum_{ps}l_{ps}x_{p1.s}+
d\sum_{ps}l_{ps}x_{p2.s} + & \\
& +e\sum_{ps}n_{ps}x_{p1.s}-\sum_{ps}n_{ps}x_{p2.s}-
k\sum_{ps}m_{ps}x_{p1.s} & = 0,
\end{array}$$
avec $A= ac+ef-ck$, $B= bc+ad+eh-kd-f$, $C= bd-h$.
Si cette relation n'est pas identiquement nulle, $a$ param\`etres
sont d\'etermin\'es. Si elle est identiquement nulle, les \'equations 
$B= C= 0$ d\'eterminent $f$ et $h$ en fonction des $a,b,c,d,e,k$.
Puis l'annulation des coefficients des $x_{pq.}$ montre que les $m_{p.}$, 
$p\neq 2$, et les $n_{q.}$, sont d\'etermin\'es par $m_{2.}$ et les $l_{r.}$, 
eux-m\^emes soumis aux relations
$$\begin{array}{rcl}
(c+ed)l_{4.} & =  & 0, \\
(c+ed)(l_{1.}+el_{2.}+l_{3.}) & =  & 0.
\end{array}$$
Enfin, la relation $B= 0$ \'equivaut \`a $(c+ed)(a-k+eb)= 0$.
Au pire, donc, $c+ed= 0$ et $7f+3$ param\`etres au moins sont 
d\'etermin\'es. On majore donc la dimension de la strate consid\'er\'ee
par le maximum de $12f+8+6a$ et $5f+5+7a$.

\medskip\noindent {\bf Cas $r=4$.}
On peut supposer que le surcro\^{\i}t de rang 
d\^{u} \`a $g$ et $g'$ se r\'epartit en $2+2$ ou $3+1$.
Dans le premier cas, on majore la dimension par $10a+8f+6$  : 
il n'y a pas de relation de sym\'etrie \`a prendre en compte.

Dans le second cas, on a des relations de d\'ependance
$$\begin{array}{rcl}
g'_{2..} & =  & ag_{1..}+bg_{2..}+cg_{3..}+dg'_{1..}+\sum_{ps}l_{ps}x_{p..s}, \\
g'_{3..} & =  & eg_{1..}+fg_{2..}+hg_{3..}+kg'_{1..}+\sum_{ps}m_{ps}x_{p..s}.
\end{array}$$
L'ouvert correspondant est donc param\'etr\'e par les $g_{pq.}$,
par les $g'_{11.}$,$g'_{14.}$, et par les coefficients $a,...,k$ et 
$l_{ps}, m_{ps}$, li\'es par la seule relation de sym\'etrie 
$g'_{23.}= g'_{32.}$. Celle-ci se traduit par l'identit\'e
$$\begin{array}{l}
Ag_{11.}+Bg_{12.}+Cg_{13.}-fg_{22.}+(b-h)g_{23.}+cg_{33.} \\
\hspace*{1.5cm} +d\sum_{ps}m_{ps}x_{p1.s}+\sum_{ps}l_{ps}x_{p3.s}-
k\sum_{ps}l_{ps}x_{p1.s}-\sum_{ps}m_{ps}x_{p2.s}= 0,
\end{array}$$
avec $A= de-ak$, $B= df-kb-e$, $C= a+dh-kc$. Si cette relation 
est identiquement nulle, les param\`etres $a,...,k$, $l_{ps}$
et $m_{ps}$ sont d\'etermin\'es par $b,d,k$ et les $m_{3s}$.
La dimension de la strate correspondante est donc major\'ee par 
le maximum de $10a+8f+8$ et $11a+f+3$.

\medskip\noindent {\bf Cas $r=5$.}
Il n'y a qu'un seul cas qui conduit au majorant $14a+4f+5$. \qed

\section{Evaluation verselle}

Cette section est consacr\'ee \`a la d\'emonstration du Th\'eor\`eme 
1.2 de l'introduction.

\begin{proposition}
On suppose $5a/2 \le b \le 4a$, $3a>11f$, $4f\le 4a-b-2$ et 
$f\ge 9a-3b$, et l'on pose $\delta = 3b-9a+f$. 
Soit $H$ un hyperplan fix\'{e} quelconque de $V$. Alors, pour $m$
en dehors d'un sous-sch\'ema de codimension $\delta+1$
de $PW$, le rang de $m_H(1)$ vaut $9a-f$.
\end{proposition}

\noindent {\em D\'{e}monstration :} 
On stratifie par le $(Z,H)$-rang $r$ de $T$,
le sch\'{e}ma $Q$ des quadruplets $(m,T,Z^{\prime },Z)$
avec $m$ dans $PW$, $Z=Z_m$, et $T$ hyperplan de $Z^{\prime}$
tel que ${\rm Im} (m_H(1))\subset T$.
On calcule la dimension des strates  $S_r$ correspondantes en
projetant sur le sch\'{e}ma des couples $(Z,T)$, ce qui donne 
 $$\dim S_r=d+\dim G_r-rb,$$
o\`u $d$ est la dimension de $PW$, et les dimensions des 
strates $G_r$ de l'espace des hyperplans de $Z'$ sont majore\'es
gr\^ace \`a la proposition 5.3. Comme $G_3$ est la strate g\'en\'erique, 
on a $\dim S_3=d-\delta-1$, et l'on s'assure pour conclure, 
gr\^ace aux estimations de 5.3, que sous les hypoth\`eses faites 
sur les entiers $a,b,f$, les autres strates n'exc\`edent
pas cette dimension. \qed

\begin{proposition}
On suppose $5a/2 \le b\le 4a$, $3a>11f$, $\delta = 3b-9a+f\ge 0$,
et de plus $4b\le 13a-5f-5$ et $5b\le 16a-2f-5$. 
Alors, pour $m$ en dehors d'un sous-sch\'{e}ma de codimension $2\delta+4$
de $PW$, le rang de $m_H(1)$ vaut au moins $9a-f-1$, et le lieu 
o\`{u} il vaut  $9a-f-1$ est lisse de codimension $\delta+1$.
\end{proposition}

\noindent {\em D\'{e}monstration :} 
On stratifie cette fois par le $(Z,H)$-rang $r$ de $T$,
le sch\'ema $Q^{\prime }$ des quadruplets $(m,T,Z^{\prime },Z)$
avec $m$ dans $PW$, $Z=Z_m$, et $T$ sous-espace de codimension deux 
de $Z^{^{\prime }}$ tel que ${\rm Im} (m_H(1))\subset T$.
On calcule la dimension des strates $S'_r$ correspondantes en
projetant sur le sch\'{e}ma des couples $(Z,T)$, ce qui donne 
 $$\dim S'_r=d+\dim G'_r-rb,$$
o\`u $d$ est \`a nouveau la dimension de $PW$, et les dimensions des 
strates $G'_r$ de la grassmannienne des sous-espaces de codimension
deux de $Z'$ sont major\'ees par la proposition 5.5. Comme $G'_6$ 
est la strate g\'en\'erique, $\dim S'_6=d-2\delta-4$, et l'on s'assure 
gr\^ace aux estimations de 5.5, que sous les hypoth\`eses faites 
sur les entiers $a,b,f$, les autres strates n'exc\`edent
pas cette dimension (cette v\'erification est fastidieuse mais 
\'el\'ementaire : on v\'erifie que la condition la plus contraignante
est donn\'ee par la strate $G'_5$, qui impose l'in\'egalit\'e
$4f\le 4a-b-5$). Ceci assure le premier point de la proposition, 
selon lequel le rang de $m_H(1)$ vaut au moins $9a-f-1$ en dehors 
d'un sous-sch\'ema de $PW$ de codimention au moins $2\delta+4$. 
 
Par ailleurs, on s'assure que sous nos hypoth\`eses, les strates 
$S_1$ et $S_2$ introduites dans la d\'emonstration de la proposition 
pr\'ec\'edente, ont une image dans $PW$ de codimension au moins 
$2\delta+4$. Reste donc, pour conclure, 
\`{a} d\'{e}montrer la lissit\'{e} du sch\'{e}ma 
$S$ des $m$ de $PW$ pour lesquels l'image de $m_H(1)$ est de 
codimension $f+1$ dans $A\otimes H.V$, 
et de $(Z_m,H)$-rang trois. Or c'est un fait 
g\'{e}n\'{e}ral concernant les
stratifications par le rang d'un sch\'{e}ma $X$, que la strate $W^i-W^{i+1}$
o\`{u} le rang du morphisme $u :k^p\rightarrow k^q$ vaut disons $s$,
s'identifie sch\'{e}matiquement \`{a} la strate correspondante $G^i-G^{i+1}$
dans $X\times G$, o\`{u} $G$ d\'{e}signe la grassmannienne des sous-espaces
de dimension $s$ de $k^q$ (cf. \cite{ACGH}, p. 69, 84). Dans le cas qui nous
occupe, il nous suffit donc de montrer que le sch\'{e}ma des couples 
$(m,T)$, avec $m$ dans $S$ et $T=Im(m_H(1))$, est lisse. Ce sch\'{e}ma
s'identifie \`{a} son tour au sch\'ema $S^{\prime }$ des quadruplets 
$(m,T,Z^{\prime },Z)$ de $Q^{\prime }$ v\'{e}rifiant les m\^{e}mes
conditions : $m$ dans $S$ et $T=Im(m_H(1))$ (et toujours $T$ de 
$(Z_m,H)$-rang $3$). Enfin, $S^{\prime }$ est lisse parce que c'est un
ouvert d'un fibr\'{e} vectoriel sur le sch\'{e}ma lisse des triplets 
$(T,Z^{\prime },Z)$ v\'{e}rifiant $T\subset Z^{\prime }=Z\cap (A\otimes
H.V)$. \qed

\medskip Rappelons que le morphisme d'\'{e}valuation 
$$H^0(E_m(1))\otimes {\cal O}\longrightarrow E_m(1)$$ 
est dit {\em versel} si les lieux o\`{u} il est de rang constant 
sont lisses de la codimension attendue (cf. \cite{EH}, 6).

\begin{theorem}\label{evvers}

On suppose $5a/2 \le b\le 4a$ et 
$$9a-3b\le f\le \min(\frac{3a-1}{11}, \frac{13a-4b-5}{5},
\frac{16a-5b-5}{2}).$$
Alors, pour $m$ g\'{e}n\'{e}rique dans $PW$, le fibr\'e vectoriel 
$E_m(1)$, noyau du morphisme correspondant de $B\otimes {\cal O}(1)$ 
vers $A\otimes {\cal O}(2)$, est
d'\'{e}valuation verselle. \end{theorem}

\noindent {\em D\'{e}monstration}  : 
On remarque d'abord que le rang du morphisme d'\'{e}valuation de 
$E_m(1)$ est \`{a} une constante pr\`{e}s la fonction $H\rightarrow
{\rm rang} (m_H(1))$ (et ce m\^{e}me au sens fort de \cite{Hi}), 
en vertu du diagramme exact  :
$$\begin{array}{ccccccccc}
&  &  &  & 0 &  & 0 &  &  \\
&  &  &  & \downarrow &  & \downarrow &  &  \\
&  &  &  & H^0(E_m(1)) & \stackrel{ev}{\longrightarrow} & E_m(1)_h &  &  \\
&  &  &  & \downarrow &  & \downarrow &  &  \\
0 & \longrightarrow & B\otimes H & \longrightarrow & B\otimes V & 
\longrightarrow & B\otimes {\cal O}_h(1) & \longrightarrow & 0 \\
&  & \hspace*{8mm}\downarrow {\scriptstyle m_H(1)} &  & 
\hspace*{7mm}\downarrow {\scriptstyle m(1)} &  & 
\hspace*{8mm}\downarrow {\scriptstyle m_h(1)} &  &  \\
0 & \longrightarrow & A\otimes H.V & \longrightarrow & A\otimes S^2V &
 \longrightarrow & A\otimes {\cal O}_h(2) & \longrightarrow & 0
\end{array}$$
On stratifie donc l'ensemble $PW\times \PP^3$ des couples 
$(m,H)$ par le rang de $m_H(1)$. D'apr\`{e}s la 
proposition pr\'ec\'edente, modulo codimension quatre au moins,
ce rang vaut $r=9a-f$, sauf le long d'un sous-sch\'{e}ma lisse de 
codimension $\delta+1$ o\`{u} il vaut $r-1$. Par le th\'{e}or\`{e}me 
de lissit\'{e} g\'{e}n\'{e}rique (en caract\'eristique nulle), 
ceci reste vrai pour $m$ g\'{e}n\'{e}rique dans $PW$. \qed

\section{Application aux courbes}

Une courbe g\'{e}n\'{e}rique $C$ de $\PP^3$ est dite {\em \`{a}
r\'{e}solution lin\'{e}aire} si son id\'{e}al $\mathcal{I}_C$ admet une
r\'{e}solution de la forme :

$$ 0\longrightarrow a. {\cal O}(-s-2)\longrightarrow b.{\cal O}(-s-1)
\longrightarrow c.{\cal O}(-s)\longrightarrow {\cal I}_C
\longrightarrow 0.$$

On a vu dans l'introduction que si le morphisme $a.\mathcal{O}
(-s-2)\rightarrow b\mathcal{O}(-s-1)$ est consid\'{e}r\'{e} comme 
morphisme $m$ de $B$ vers $A\otimes V,$ alors $m(1)$ est de rang 
au plus $3b+a-1$. En particulier, si $b\leq 3a$, le morphisme $m$ 
ne peut \^{e}tre g\'{e}n\'{e}rique.

\begin{definition}

Nous dirons qu'une courbe g\'{e}n\'{e}rique \`{a} r\'{e}solution
lin\'{e}aire de type $(a,b)$ est {\em pr\'{e}dominante} (rlp) 
si l'image de $m(1)$ est un sous-espace g\'{e}n\'{e}rique 
(en ce sens que son image dans le
quotient de la grassmannienne par $Gl(A)$ est le point g\'{e}n\'{e}rique) de
dimension $\min (10a,3b+a-1)$ de $A\otimes S^2V$.

\end{definition}

Rappelons que les courbes g\'{e}n\'{e}riques \`{a}
r\'{e}solution lin\'{e}aire dominante introduites dans \cite{EH} 
pour $b>3a$, ont $m(1)$ surjectif et sont donc  rlp.

\begin{theorem}
Pour $b>(3-\frac{1}{11})a+\frac{9}{11}$ et $a\ge 7$, 
il existe dans $\PP^3$ une unique courbe g\'{e}n\'{e}rique lisse, 
g\'{e}om\'{e}triquement connexe, \`a r\'esolution lin\'eaire 
pr\'edominante de type $(a,b)$.
\end{theorem}

\noindent {\em D\'{e}monstration :} 
Le cas $b>3a$ a \'{e}t\'{e} trait\'{e} dans \cite{EH}. 
On suppose donc $b\le 3a$, et l'on pose $c=b-a+1$, $s=b-2a$ et
$f=9a-3b+1$.  On v\'{e}rifie
sans peine que sous l'hypoth\`ese faite sur $a$ et $b$, 
on peut appliquer \ref{evvers}. Consid\'{e}rons donc la matrice $m$ et le
fibr\'{e} $E_m(1)$ donn\'{e}s par \ref{evvers}. On sait que $h^0(E_m(1))=c$
et que $E_m(1)$ est d'\'{e}valuation verselle. En dualisant le morphisme
d'\'{e}valuation on obtient un morphisme versel de $E_m(1)^{\vee}$ vers 
$c.{\cal O}$. Comme $E_m(1)^{\vee}$ est de rang $c-1$, le conoyau de
ce morphisme
est de la forme ${\cal I}_C(s)$ o\`{u} ${\cal I}_C$ est l'id\'{e}al
d'une courbe lisse $C,$ rlp de type $(a,b)$. 

Cette courbe est g\'{e}om\'{e}triquement 
connexe puisque $h^1({\cal I}_C)=0$. En effet,
$$H^1({\cal I}_C)=H^2(E_m(1)^{\vee}(-s)),$$ 
qui est le noyau du morphisme induit par $m$ de $a.H^3({\cal O}(-s-2))$
vers  $b.H^3({\cal O}(-s-1))$. Pour montrer que ce morphisme est injectif,
il suffit de s'assurer que son dual de Serre, qui est le morphisme 
induit par $m$  de $b.H^0({\cal O}(s-3))$ vers $a.H^0({\cal O}(s-2))$,
 est surjectif. Mais cela r\'{e}sulte du point $2$ du th\'eor\`eme 4.2.

Observons enfin que $C$ est g\'{e}n\'{e}rique. En effet, une
g\'{e}n\'{e}ralisation $C^{\prime }$ de $C$ a une r\'{e}solution du
m\^{e}me type (cf. \cite{EH}, d\'{e}monstration de 7.1). Soit $m^{\prime }$
la matrice correspondante. L'image de $m^{\prime }(1)$ est encore de
codimension $f$ (car $h^0(E_{m^{\prime }}(1))=c$), et g\'{e}n\'{e}ralise par
cons\'{e}quent celle de $m(1)$. Cette image est donc le sous-espace
g\'{e}n\'{e}rique de codimension $f$ de $A\otimes S^2V$  :  par
cons\'{e}quent, $C^{\prime }$ est rlp. Or toute courbe $C^{\prime }$
g\'{e}n\'{e}rique rlp de type $(a,b)$ est \'{e}gale \`{a} $C.$ En effet, 
$m^{\prime }=m$ \`{a} extension des scalaires pr\`{e}s, puisque $m$ est la
matrice g\'{e}n\'{e}rique telle que $m(1)$ ait pour image le sous-espace
g\'{e}n\'{e}rique de codimension $f$ de $A\otimes S^2V.$ Comme $C^{\prime }$
est le lieu de d\'{e}g\'{e}n\'{e}rescence du morphisme d'\'{e}valuation de 
$E_{m^{\prime }}(1),$ elle est bien \'{e}gale \`{a} $C$. \qed

\medskip\noindent {\bf Remarque 7.3} L'\'enonc\'e pr\'ec\'edent est
tout pr\`es d'\^etre optimal. Supposons en effet que 
$b<(3-\frac{1}{11})a+\frac{1}{3}$. Alors d'apr\`es la remarque 5.4, si
$H$ est fix\'e et $Z$ g\'en\'erique de codimension $f=9a-3b+1$ dans 
$A\otimes S^2V$, l'image de $m_H(1)$ est incluse dans un hyperplan de 
$Z^{\prime}$. Ceci exclut l'existence d'une courbe g\'en\'erique 
rlp de type $(a,b)$.

\end{document}